\UseRawInputEncoding
\documentclass[a4paper]{amsart}
\usepackage{amscd,amsthm,amsfonts,amssymb,amsmath,supertabular}
\usepackage{mathrsfs}
\usepackage[]{hyperref}
\usepackage{fullpage}
\usepackage[all]{xy}

    \newcommand{\BA}{{\mathbb {A}}} 
    \newcommand{\BC}{{\mathbb {C}}} 
     \newcommand{\BF}{{\mathbb {F}}}
    \newcommand{\BG}{{\mathbb {G}}}

     \newcommand{\BN}{{\mathbb {N}}}
     \newcommand{\BP}{{\mathbb {P}}}
    \newcommand{\BQ}{{\mathbb {Q}}} \newcommand{\BR}{{\mathbb {R}}}

     \newcommand{\BZ}{{\mathbb {Z}}}

    \newcommand{\CO}{{\mathcal {O}}}

     \newcommand{\RH}{{\mathrm {H}}}

    \newcommand{\fa}{{\mathfrak{a}}} \newcommand{\fb}{{\mathfrak{b}}}

     \newcommand{\fp}{{\mathfrak{p}}}

     \newcommand{\fD}{{\mathfrak{D}}}

     \newcommand{\fP}{{\mathfrak{P}}}

    \newcommand{\Cl}{{\mathrm{Cl}}}

    \newcommand{\Div}{{\mathrm{Div}}} \renewcommand{\div}{{\mathrm{div}}}
    \newcommand{\End}{{\mathrm{End}}} 
    
    \newcommand{\Fr}{{\mathrm{Fr}}}

    \newcommand{\Gal}{{\mathrm{Gal}}} \newcommand{\GL}{{\mathrm{GL}}}
    
    \newcommand{\Hom}{{\mathrm{Hom}}}
    
    \newcommand{\id}{{\mathrm{id}}}\renewcommand{\Im}{{\mathrm{Im}}}

    \newcommand{\Jac}{{\mathrm{Jac}}}\newcommand{\Ker}{{\mathrm{Ker}}}

    \newcommand{\loc}{{\mathrm{loc}}}

    \newcommand{\ord}{{\mathrm{ord}}} \newcommand{\rk}{{\mathrm{rank}}}
     
    \newcommand{\pr}{{\mathrm{pr}}} 
    
    \renewcommand{\mod}{\ \mathrm{mod}\ }\renewcommand{\Re}{{\mathrm{Re}}}

    \newcommand{\Sel}{{\mathrm{Sel}}}

    \newcommand{\sgn}{{\mathrm{sgn}}}

    \newcommand{\tr}{{\mathrm{tr}}}\newcommand{\tor}{{\mathrm{tor}}}
    
    \newcommand{\ur}{{\mathrm{ur}}}

        \newcommand{\Tr}{\mathrm{Tr}}

    \font\cyr=wncyr10

    \newcommand{\Sha}{\hbox{\cyr X}}
    \newcommand{\wh}{\widehat}

    \newcommand{\ov}{\overline}

    \newcommand{\sk}{\medskip}
    \newcommand{\lra}{\longrightarrow}
    \newcommand{\ra}{\rightarrow}

    \newcommand{\s}{\sk\noindent}
    
    \newcommand{\N}{\mathrm{N}}

    \newcommand{\m}{\mathrm{m}}                     
                            
                              \newcommand{\hilbert}[4]{\left(\frac{#1,#2}{#3;#4}\right)}

    \theoremstyle{plain}
    \newtheorem{thm}{Theorem}[section] \newtheorem{coro}[thm]{Corollary}
    \newtheorem{lem}[thm]{Lemma}  \newtheorem{prop}[thm]{Proposition}

\theoremstyle{remark} \newtheorem{remark}{Remark}[section]
\theoremstyle{remark} 
\theoremstyle{remark} 

    \newcommand{\et}{\mathrm{\acute{e}t}}

    \numberwithin{equation}{section}

\begin{document}
\title{Root numbers and Selmer groups for the Jacobian varieties of Fermat curves}
\author{Jie Shu}

\address{School of Mathematical Sciences, Tongji University, Shanghai 200092}

\subjclass[2010]{Primary 11G40; 14G05; 14G10}
\thanks{The author is supported by NSFC-11701092.}

\email{shujie@tongji.edu.cn}\maketitle

\begin{abstract}
Let $p$ be an odd prime number. Let $K$ be the $p$-th cyclotomic field and $F$ its maximal real subfield. We give general formulae of the root numbers of the Jacobian varieties of the Fermat curves $X^p+Y^p=\delta$ where $\delta$ is an integer. As an application of these general formulae, we derive the equidistribution of the root numbers for the families of Jacobian varieties of the Fermat curves. When $p\nmid \delta$, we bound the Selmer groups of these Jacobian varieties. Moreover, if $p$ is regular and all prime ideals of $K$ dividing $\delta$ are inert in $K/F$, the Selmer groups are explicitly determined and we verify the $p$-parity conjectures of these Jacobian varieties. We also give an asymptotic lower bound for the number of  Fermat Jabobians for which the $p$-parity conjecture holds.
\end{abstract}

 \s{\bf Notations.} Let $\mu_p$ be the group of $p$-th root of unity with a generator $\omega$. Let $K=\BQ(\omega)$ be the $p$-th cyclotomic field and $F=\BQ(\omega+\omega^{-1})$ the maximal real subfield of $K$. Let $\CO_K$ resp. $\CO$ be the ring of integers of the field $K$ resp. $F$. Let $\BA$ be the ad\'ele ring of $F$ and $\BA_f$ its finite part. For any $\BZ$-module $M$, we denote by $\wh{M}=M\otimes_\BZ \wh{\BZ}$ and $\wh{\BZ}=\prod_p \BZ_p$. For example, $\wh{F}=\BA_f$. Let $|\cdot|_\BA: \BA^\times \lra \BR_+^\times$ denote the standard ad\'elic absolute value so that $d(ab)=|a|_\BA db$ for any Haar measure $db$ on $\BA$. Let $|\cdot |_v$ denote the absolute value on $F_v^\times$ for each place $v$ of $F$ such that $|x|_\BA=\prod_v |x_v|_v$ for any $x=(x_v)\in \BA^\times$. For a finite place $v$, sometimes we also denote by $v$ its corresponding prime ideal and  $q_v=\# \CO/v$.  For $x\in F^\times$, denote by  $\ord_v(x)$ the additive valuation of $x$ at $v$ so that $\ord_v(v)=1$. We denote by $\infty$  the set of infinite places of $F$.  Let $\BA_K=\BA\otimes_F K$ be the ad\'ele ring of $K$. The notations $|\cdot |_{\BA_K}, |\cdot |_V,q_V,$ and $\ord_V(\cdot)$ are similarly defined for $K$ where $V$ is a finite place of $K$.

\section{Introduction}
Let $p$ be an odd prime number. Let $\mu_p$ be the group of $p$-th root of unity, and let $\omega$ be a primitive $p$-th root of unity. Let $K=\BQ(\mu_p)$ be the $p$-th cyclotomic field. Let $\delta\in \BQ^\times/\BQ^{\times p}$ be an integer and let $F_\delta$ be the Fermat curve defined by
\[X^p+Y^p= \delta.\]
We have two algebraic automorphisms  $\tau_1$ and $\tau_2$  of $F_\delta$ as follows:
\[\tau_1(X,Y)\mapsto(\omega X, Y),\quad \tau_2(X,Y)\mapsto( X, \omega Y).\]

For integers $r,s,t >0$ satisfying $r+s+t=p$, let $C_{r,s,t;\delta}$ be the quotient of $F_\delta$ by the group $\langle \tau_1^s\cdot\tau_2^{-r}\rangle $. Then $C_{r,s,t;\delta}$ is the projective curve of genus $g=(p-1)/2$, which is  defined by the affine equation
\[y^p=x^r( \delta-x)^s,\]
and, in the projective coordinates, the quotient map is given by
\[\varphi_{r,s,t}:F_\delta\lra C_{r,s,t;\delta},\quad (X,Y,Z)\mapsto (X^p, X^rY^sZ^t,Z^p).\]
Denote by $J_{r,s,t;\delta}$ the Jacobian variety of the curve $C_{r,s,t;\delta}$. We have a complex multiplication of $\mu_p$ on $C_{r,s,t;\delta}$:
\[[\omega](x,y)=(x,\omega y).\]
This induces a complex multiplication $[\ ]:\CO_K\hookrightarrow\End(J_{r,s,t; \delta})$. 

Let $\varphi$ be the product morphism
\[\varphi=(\varphi_{r,1,p-r-1})_{r=1,2,\cdots,p-2}:F_\delta\lra \prod_{r=1}^{p-2}C_{r,1,p-r-1;\delta}.\]
It is known from \cite{Fad61,Fad61b,Rohrlich78} that the push-forward map 
\[\varphi_*:\Jac(F_\delta)\ra \prod_{r=1}^{p-2} J_{r,1,p-r-1;\delta}\]
gives an isogeny of abelian varieties over $\BQ$.  This isogeny gives a rough decomposition of the Jacobian variety $\Jac(F_\delta)$ and the abelian varieties $J_{r,s,t;\delta}$ will be referred as the Jacobian varieties of the Fermat curves.

The curve 
\[C_{r,s,t;\delta}:y^p=x^r(\delta-x)^s\]
is defined over $\BQ$.  Let  $\Lambda(s,{J_{r,s,t;\delta}}_{/\BQ})$ be the complete  global $L$-function
of the Jacobian variety ${J_{r,s,t;\delta}}_{/\BQ}$. Weil \cite{Weil52} proved this $L$-function is modular and can be  identified with the $L$-function of a Hecke character over $K$. Therefore $L(s,{J_{r,s,t;\delta}}_{/\BQ})$ is entire and satisfies Hecke's functional equation 
\[\Lambda(s,{J_{r,s,t;\delta}}_{/\BQ})=\epsilon_{r,s,t;\delta}\Lambda(2-s,{J_{r,s,t;\delta}}_{/\BQ})\]
where $\epsilon_{r,s,t;\delta}=\pm 1$ is the root number (or epsilon factor). Combined with the Taylor expansion at $s=1$, the functional equation implies 
\[\epsilon_{r,s,t;\delta}=(-1)^{\ord_{s=1 }\Lambda(s,{J_{r,s,t;\delta}}_{/\BQ})}.\]

By the theorem of Mordell-Weil, 
\[J_{r,s,t;\delta}(\BQ)\simeq \BZ^{R_{r,s,t;\delta}}\times J_{r,s,t;\delta}(\BQ)_\tor,\]
where $R_{r,s,t;\delta}\geq 0$ is an integer. Then the Birch and Swinnerton-Dyer conjecture predicts that the sign of the epsilon factor $\epsilon_{r,s,t;\delta}$ is compatible with the parity of the Mordell-Weil rank $R_{r,s,t;\delta}$, i.e.
\[(-1)^{R_{r,s,t;\delta}}=\epsilon_{r,s,t;\delta}. \]

The field extension $K/\BQ$ is totally ramified at $p$ and let $\Pi=\omega-\ov{\omega}$ be a primitive element of $K$ above $p$. Denote by $\Sel\left({J_{r,s,t;\delta}}_{/K},\Pi\right)$ the Selmer group associated to the $K$-morphism $[\Pi]:J_{r,s,t;\delta}\ra J_{r,s,t;\delta}$. Then $\Sel\left({J_{r,s,t;\delta}}_{/K},\Pi\right)$ is a finite dimensional vector space over $\BF_p=\CO_K/(\Pi)$. Denote  
\[S_{r,s,t;\delta}= \dim_K {J_{r,s,t;\delta}(K)\otimes_{\CO_K}K}+\rk_{\BF_p}\Sha\left({J_{r,s,t;\delta}}_{/K}\right)[\Pi],\]
where $\Sha\left({J_{r,s,t;\delta}}_{/K}\right)$ denotes the Shafarevich-Tate group of the abelian variety ${J_{r,s,t;\delta}}_{/K}$.
Then $S_{r,s,t;\delta}$ is essentially the dimension of the Selmer group. Precisely, it is known from Proposition \ref{torsion} and Corollary \ref{selmer-rank} that
\[S_{r,s,t;\delta}=\dim_{\BF_p}\Sel\left({J_{r,s,t;\delta}}_{/K},\Pi\right)-1.\]
As one easily show that
\[\dim_{\BQ}J_{r,s,t;\delta}(\BQ)\otimes_\BZ\BQ=\dim_K J_{r,s,t;\delta}(K)\otimes_{\CO_K}K,\]
and there is conjecturally a non-degenerate alternating pairing on $\Sha\left({J_{r,s,t;\delta}}_{/K}\right)$, one should have
\[R_{r,s,t;\delta}\equiv S_{r,s,t;\delta}\mod 2.\]
Hence by the conjectures of Birch and Swinnerton-Dyer, we should have the parity:
\begin{equation}\label{pc}
\epsilon_{r,s,t;\delta}=(-1)^{S_{r,s,t;\delta}}.
\end{equation}

When $p=3$ and $r=s=t=1$, both $F_\delta$ and $C_{1,1,1;\delta}$ are elliptic curves. The parity conjecture (\ref{pc}) in this case is extensively verified by the works of Birch, Stephens, Liverance and  Satg\'e \cite{BS66, Stephens68, Liver95,Satge}. Indeed, Birch, Stephens and Liverance \cite{BS66, Stephens68, Liver95} calculated the general formulae of the root numbers of these elliptic curves and Satg\'e \cite{Satge} calculated the Selmer groups of these elliptic curves by descent method via $3$-isogenies. For the cases of  $\delta=1$ and regular primes $p\geq 5$ , the parity conjecture $(\ref{pc})$ is verified, at least when $r=s=1$,  by the works of Gross, Rohrlich and Faddeev \cite{Gro-Roh78, Rohr92, Fad61}. Gross and Rohrlich \cite{Gro-Roh78, Rohr92} calculated the root numbers of the Jacobian varieties of the Fermat curves $X^p+Y^p=1$ for odd primes $p$, while Faddeev \cite{Fad61} calculated the Selmer groups of these Jacobian varieties when $p\geq 5$ is regular. 

For a $p$-th power-free integer $\delta$, we denote the following Hypothesis ($\mathrm{\bf{I}}$) for $\delta$:
\[\label{I}\text{any prime ideal of $K$ dividing $\delta$ is inert in $K/F$.}\leqno{\mathrm{(\bf{I})}}\]
We shall prove the following 
\begin{thm}\label{parity-intro}
Suppose $p\geq 5$ is a regular prime, i.e.  $\Cl(K)[p]=0$, and $(r,s,t)$ is a triple of integers with $r,s,t>0$ and $r+s+t=p$. Suppose $\delta$ is a $p$-th power-free integer satisfying the Hypothesis $(\mathrm{\bf{I}})$. Then \[\epsilon_{r,s,t;\delta}=(-1)^{S_{r,s,t;\delta}}.\]
\end{thm}

The above theorem is Theorem \ref{parity} where the root number is given in terms of the root number $\epsilon(\pi_{r,s,t;\delta})$ of the automorphic representation $\pi_{r,s,t;\delta}$ associated to the Jacobian variety $J_{r,s,t;\delta}$ via the Eichler-Shimura construction, and it is known that $\epsilon_{r,s,t;\delta}=\epsilon(\pi_{r,s,t;\delta})$. In order to prove this theorem, we shall calculate the root numbers and bound the Selmer groups respectively.

For a fixed odd prime $p$, we know from Proposition \ref{density} that primes satisfying $(\mathrm{\bf{I}})$ have natural density $1-2^{-n}\geq 1/2$ among all  primes where $n=\ord_2(p-1)$. This density converges to one rapidly. By an variant \cite{Delang56} of Ikehara's Tauberian theorem, one can give the  asymptotic approximation of the number of $p$-th power-free integers $\delta$ satisfying $(\mathrm{\bf{I}})$. When $p$ is regular, by Theorem \ref{parity-intro}, this also gives the asymptotic lower bound for the number of  Fermat Jabobians for which the parity conjecture  (\ref{pc}) holds.

For the prime $p$  and the triple $(r,s,t)$ with $r,s,t>0$ and $r+s+t=p$, we define
\[N_{r,s,t;\delta}(p,X)=\#\{\delta < X : |\delta|\leq X, \epsilon_{r,s,t;\delta}=(-1)^{S_{r,s,t;\delta}}\}.\]
One has the following asymptotic lower bound for $N_{r,s,t;\delta}(p,X)$.

\begin{thm}
Suppose $p\geq 5$ is a regular prime and $(r,s,t)$ is a triple of integers with $r,s,t>0$ and $r+s+t=p$. As $X\ra \infty$, we have
\[N_{r,s,t;\delta}(p,X)\gg \frac{X}{\log^{2^{-n}}X},\]
where $n=\ord_2(p-1)$.
\end{thm}

\begin{remark}
It is conjectured that there are infinitely many regular primes. More precisely Siegel (\cite{Siegel64}) conjectured that $e^{-1/2}$, or about $60.65\%$, of all prime numbers are regular (see \cite{HHO}, \cite{Gras17} for numerical evidence). It seems quite possible, based on numerical evidence, that the valuation $n=\ord_2(p-1)$ can be arbitrarily large for regular primes $p$.  Indeed, for any given $n$, by Chebotarev's density theorem, there is a positive proportion of primes $p\equiv 1\mod 2^n$. If the regular primes are equidistributed, there should be many  regular primes (with density $>1/2$)  among primes $p\equiv 1\mod 2^n$.   The  following table is communicated by Prof. G. Gras. 


\begin{center}
\tablefirsthead{
\multicolumn{4}{c}{{\bf Table.} Regular primes $p$ with increasing $n$.}\\
\multicolumn{4}{c}{}\\
\hline  $p$ & $n$ & density of primes satisfying (I) & lower bound for $N_{r,s,t;\delta}(p,X)$ \\
\hline }

\tablehead{
\hline  $p$ & $n$ & $d$ & lower bound \\
\hline }

\tabletail{\hline}

\tablelasttail{\hline}
\begin{supertabular}{|c|c|c|c|}
$5$&$2$& $0.7500$&$X/\log^{2^{-2}}X$\\
$17$&$4$& $0.9375$&$X/\log^{2^{-4}}X$\\
$769$&$8$& $0.9960$&$X/\log^{2^{-8}}X$\\
$12289$&$12$& $0.9998$&$X/\log^{2^{-12}}X$\\
$40961$&$13$& $0.9999$&$X/\log^{2^{-13}}X$\\
\hline
\end{supertabular}
\end{center}

\end{remark}

In \S 2 we sketch the basics on Jacobi sums and algebraic Hecke characters and L-functions associated to the Jacobian varieties of Fermat curves. We also roughly determine the $p$-torsion groups of the Mordell-Weil groups of the Jacobian varieties over the cyclotomic field $K$.

In \S 3 we calculate the general formulae for the root numbers $\epsilon(\pi_{r,s,t;\delta})$ and the results are summarized in Theorem \ref{root-number}. The formulae given here are unconditional, that is, we do not require the hypothesis for $\delta$ as in Theorem \ref{parity-intro}. These results are compatible with the works of Birch, Stephens and Liverance \cite{BS66, Stephens68, Liver95} for the case $p=3$ and Gross-Rohrlich \cite{Gro-Roh78, Rohr92} for the case $\delta=1$. 

As an application of these general formulae,  following the strategy and results in \cite{Mai93, Tomasz08}, we obtain the following density of root numbers of the families $J_{r,s,t;\delta}$ with $\delta$ varying in the set of  $p$-th power-free integers.  Let $\BN^{(p)}$ be the set of $p$-th power-free integers.

\begin{thm}
\[\lim \limits_{X\ra +\infty} \frac{\#\{\delta\in \BN^{(p)}:\,\delta\leq X\textrm{ and } \epsilon_{r,s,t;\delta}=+1\}}{\#\{\delta\in \BN^{(p)}: \,\delta\leq X \}}=\frac{1}{2}.\]
\end{thm}

In \S 4 we bound the Selmer groups $\Sel\left({J_{r,s,t;\delta}}_{/K},\Pi\right)$.  The Selmer group $\Sel\left({J_{r,s,t;\delta}}_{/K},\Pi\right)$ consists of 1-cocycles in $\RH^1(K,J_{r,s,t;\delta}[\Pi])$ which  lie in the image of the local Kummer maps
\[\kappa_V:J_{r,s,t;\delta}(K_V)/\Pi J_{r,s,t;\delta}(K_V)\hookrightarrow \RH^1(K_V,J_{r,s,t;\delta}[\Pi])\] 
for each place $V$ of $K$. In order to bound the Selmer gorups, we compute the images of the local Kummer  maps via Weil's point of view of the Kummer maps. This computation of local Kummer images is relatively easy to carry out for  the places not above $p$ (Theorem \ref{desc-off-p}). But for the place above $p$, it is complicated and we generalize the arguments in Faddeev \cite{Fad61} (Theorem \ref{desc-at-p1}).

Let $R=\CO_K\left[\frac{1}{p\delta}\right]$ and then $R^\times$ is the groups of $p\delta$-units in $K$. We have an exact sequence
\begin{equation}\label{class}
0\ra \Ker(\alpha)\ra \Sel\left({J_{r,s,t;\delta}}_{/K}, \Pi\right)\xrightarrow{\alpha} \Cl(R)[p],
\end{equation}
where
\[\Ker(\alpha)=\{x\in R^\times/R^{\times p}: x_V\in \Im(\delta_V) \textrm{ for all } V\}.\]
Hence, modulo the class group $\Cl(R)[p]$, which is controlled by $\Cl(K)[p]$, the investigation of the Selmer groups relies heavily on the $\BF_p$-space $\Ker(\alpha)$.   

Let $i(p)$ denote the irregularity of $p$, i.e. the number of Bernoulli numbers $B_{i}$, with $i$ even and $2\leq i\leq p-3$, such that $\ord_p(B_{i})>0$. 
A similar investigation of the space $\Ker(\alpha)$  as in Mccallum \cite{McCallum92} enables us to bound the Selmer groups using the number $k(\delta)$ of distinct prime ideals of $K$ dividing $\delta$,  the irregularity $i(p)$  and the class group $\Cl(K)[p]$ (Theorem \ref{selmer-bound}) as follows.

\begin{thm}\label{selmer-bound-intro}
Suppose $p\geq 5$ is a prime, $(r,s,t)$ is a triple of integers with $r,s,t>0$ and $r+s+t=p$, and $\delta$ is a $p$-th power-free integer such that $p\nmid \delta$. Then
\[\dim_{\BF_p}\Sel\left({J_{r,s,t;\delta}}_{/K}, \Pi\right)\leq k(\delta)+\frac{p-3}{4}+i(p)+\dim_{\BF_p}\Cl(K)[p].\]
\end{thm}

By Ribet's converse to Herbrand's theorem \cite{Ribet76},
\[i(p)\leq \dim_{\BF_p}\Cl(K)[p].\]
Therefore, we have
\begin{coro}
Let the notations be as in Theorem \ref{selmer-bound-intro}.  Then
\[\dim_{\BF_p}\Sel\left({J_{r,s,t;\delta}}_{/K}, \Pi\right)\leq k(\delta)+\frac{p-3}{4}+2\dim_{\BF_p}\Cl(K)[p].\]
\end{coro}

If the prime $p$ is regular, i.e. $\Cl(K)[p]=0$, then $\Cl(R)[p]=0$. By the exact sequence (\ref{class}), the Selmer groups are exactly the spaces $\Ker(\alpha)$, and the structure of the Selmer groups can be explicitly determined.

\begin{thm}\label{Selmer-intro}
Suppose $p\geq 5$ is a regular prime, the triple $(r,s,t)=(r,1,p-r-1)$ with $1\leq r\leq p-2$ and $\delta$ is a $p$-th power-free integer satisfying the Hypothesis $(\mathrm{\bf{I}})$. Then
\[\dim_{\BF_p}\Sel\left({J_{r,s,t;\delta}}_{/K},\Pi\right)=k(\delta)+\left\{\begin{aligned}
\frac{1}{4}\left(p-3+2\left(\frac{B}{p}\right)\right),\quad \textrm{ if } p\equiv 1\mod 4;\\
\frac{1}{4}\left(p-5-2\left(\frac{B}{p}\right)\right),\quad \textrm{ if } p\equiv 3\mod 4.
\end{aligned}\right.\]

\end{thm}
These dimension formulae are established in Theorem \ref{Selmer} and the notaion $B$  is as  that there. Moreover, the explicit structure of the Selmer groups  are given in terms of generators of $R^\times/R^{\times p}$ there. These dimension formulae  are compatible with those in \cite[Proposition 4.1]{Gro-Roh78} for the case $\delta=1$ as expected (although there is a sign error in the second formula of \cite[Proposition 4.1]{Gro-Roh78} ).

In the final section, we note that under the condition of Theorem \ref{parity-intro}, the explicit structures of the Selmer group $\Sel\left({J_{r,s,t;\delta}}_{/K},\Pi\right)$ are determined.  Comparing the dimension formulae of Selmer groups in Theorem \ref{Selmer-intro} with the formulae for root numbers in Theorem \ref{root-number} and noting Corollary \ref{selmer-rank} that
\[S_{r,s,t;\delta}=\dim_{\BF_p}\Sel\left({J_{r,s,t;\delta}}_{/K}, \Pi\right)-1,\] 
we verify the parity conjecture Theorem \ref{parity-intro}. By an variant \cite{Delang56} of Ikehara's Tauberian theorem, we also give the asymptotic lower bound for the number of  Fermat Jabobians for which the parity conjecture  (\ref{pc}) holds.

\noindent {\bf Acknowledgements.} The author would like to thank Professor Georges Gras, Benedict Gross for helpful communications.

\section{The Jacobian Varieties of Fermat Curves}
\subsection{Jacobi sums and algebraic Hecke characters}
The readers can refer to \cite{Weil52, CM88} for the basics of the Jacobi sums. Let $r,s,t >0$ be integers satisfying $r+s+t=p$. Let $V$ be any finite place of $K$ not dividing $p$.
The p-th power residue character $\chi_V:(\CO_K/V)^\times\ra \mu_p$ modulo $V$  is given by
\[\chi_V(x)\equiv x^{\frac{q_V-1}{p}}\mod V.\]
The Jacobi sum
\[j_{r,s,t}(V)=-\sum_{\begin{subarray}{c}x\in \CO_K/V\\x\neq 0,1\end{subarray}}\chi_V(x)^r\chi_V(1-x)^s\]
is an integer in $\CO_K$. The set
\[\Phi_{r,s,t}=\left\{h\in (\BZ/p\BZ)^\times|\, \left\langle \frac{hr}{p}\right\rangle+ \left\langle \frac{hs}{p}\right\rangle+\left\langle \frac{ht}{p}\right\rangle=1\right\}\]
 can be identified with the CM-type of $J_{r,s,t;1}$ through the Artin isomorpshim
\[ \sigma:(\BZ/p\BZ)^\times\simeq \Gal(K/\BQ),\quad \sigma_h(\omega)=\omega^h.\]
If we put 
\[\Phi_{r,s,t}^{-1}=\{h^{-1}\in (\BZ/p\BZ)^\times| \, h\in \Phi_{r,s,t}\},\]
then the CM type $(K,\Phi_{r,s,t}^{-1})$ is the reflex of the CM type $(K,\Phi_{r,s,t})$. 
The Stickelberger relation holds
\begin{equation}\label{Stic}
(j_{r,s,t}(V))=\prod_{h\in \Phi_{r,s,t}}V^{\sigma_h^{-1}},
\end{equation}
and it follows that $j_{r,s,t}(V)$ has absolute value $q_V^{1/2}$ for any complex embedding of $K$ in $\BC$.

We extend the definition of $j_{r,s,t}(V)$ to all ideals prime to $p$ in $K$ by the condition
\[j_{r,s,t}(\fa\fb)=j_{r,s,t}(\fa)j_{r,s,t}(\fb).\]
Weil proved in \cite{Weil52}  that the resulting function $j_{r,s,t}(\fa)$ is indeed an algebraic Hecke character on $K$ with values in $K$ of conductor dividing $\Pi^2$ and infinity type $\Phi^{-1}_{r,s,t}$, i.e. for any $x\in K^\times $ such that $x\equiv 1\mod \Pi^2$,
\[j_{r,s,t}((x))=x^{\Phi_{r,s,t}^{-1}}=\prod_{\sigma\in \Phi_{r,s,t}^{-1}}x^{\sigma}.\]
Here $\Pi=\omega-\ov{\omega}$ is a primitive element of $K$ at the place above $p$.
The infinity type
\[\Phi_{r,s,t}^{-1}:K^\times\lra K^\times,\quad x\mapsto \prod_{\sigma\in \Phi_{r,s,t}^{-1}} x^{\sigma},\]
induces a homomorphism, by linearity, on id\'eles:
\[\Phi_{r,s,t}^{-1}:\BA_K^\times\lra \BA_K^\times.\]
Let $\tau_0\in \Hom(K,\BC)$ be a fixed element and we identify $K$ as a subfield of $\BC$ via $\tau_0$.  Let 
\[\Phi_{r,s,t;\tau_0}^{-1}:\BA_K^\times\lra \BA_K^\times\xrightarrow{\tau_0}\BC\]
be the composition of $\Phi_{r,s,t}^{-1}$ with the projection through $\tau_0$.
It is well known \cite{Weil52} that 
\[\phi_{r,s,t}=j_{r,s,t}\cdot(\Phi_{r,s,t;\tau_0}^{-1})^{-1}:K^\times\backslash\BA_K^\times\ra \BC^\times\]
is the Hecke character associated to the CM abelian variety ${J_{r,s,t;1}}_{/K}$ in the sense of \cite[Proposition 19.9]{Shimura-CM}. Generally, if $\delta\in K^\times/K^{\times p}$, then according to Weil \cite{Weil52} the abelian variety ${J_{r,s,t;\delta}}_{/K}$ also has CM type $(K,\Phi_{r,s,t})$ and the associated Hecke character  is $\phi_{r,s,t;\delta}=\chi_ \delta^{r+s}\cdot j_{r,s,t}\cdot(\Phi_{r,s,t;\tau_0}^{-1})^{-1}$ where $\chi_ \delta$ is the character on ideals prime to $ \delta$ such that $\chi_ \delta(V)=\chi_V( \delta)$. It follows from the Stickelberger relation (\ref{Stic}) that the Hecke characters $$\phi_{r,s,t; \delta}=(\chi_ \delta^{r+s}j_{r,s,t}\cdot(-\Phi^{-1}_{r,s,t}))_{\tau_0}|\cdot |_{\BA_K}^{1/2}:K^\times\backslash \BA_K^\times\lra \BC^\times$$ are all unitary.

\subsection{The automorphic representations and L-functions}\label{L-function}
Fix the triple $(r,s,t)$ with $r,s,t>0$ and $r+s+t=p$. Suppose $\delta \in F^\times/F^{\times p}$. Then both $C_{r,s,t;\delta}$ and $J_{r,s,t;\delta}$ are defined over $F$. Let $\pi_{r,s,t;\delta}$ be the irreducible automorphic representation associated to ${J_{r,s,t;\delta}}_{/F}$ by the Eichler-Shimura theory. Explicitly $\pi_{r,s,t;\delta}$ is the representation on $\GL_2(\BA)$ locally constructed from $\phi_{r,s,t;\delta}$ using Weil representation. 
The L-function $L(s,\pi_{r,s,t;\delta})$ is entire and satisfies a functional equation
\[L(s,\pi_{r,s,t;\delta})=\epsilon(\pi_{r,s,t;\delta})L(1-s,\pi_{r,s,t;\delta}),\]
where the epsilon factor $\epsilon(\pi_{r,s,t;\delta})=\pm 1$.
By \cite[Theorem 4.7]{JL70}, we know
\[L(s,\pi_{r,s,t;\delta})=L(s,\phi_{r,s,t;\delta}).\]

 Furthermore, suppose $\delta\in \BQ^\times/\BQ^{\times p}$ is an integer. The curve 
\[C_{r,s,t;\delta}:y^p=x^r(\delta-x)^s\]
is defined over $\BQ$ and has good reduction at all rational primes $\ell\nmid p\delta$. Let 
\[Z({C_{r,s,t;\delta}}_{/\BF_\ell},T)=\frac{P_\ell(T)}{(1-T)(1-\ell T)}\]
be the Zeta function of the reduced curve ${C_{r,s,t;\delta}}_{/\BF_\ell}$. Then $P_\ell(T)$ is a polynomial with integral cofficients of degree $2g=p-1$.

The Jacobian variety ${J_{r,s,t;\delta}}_{/\BQ}$ has a  global $L$-function
\[L(s,{J_{r,s,t;\delta}}_{/\BQ})=\prod_{\ell\nmid p\delta}P_\ell(\ell^{-s})^{-1}.\]
This Euler product converges for $\Re(s)>3/2$. Weil \cite{Weil52} proved in this region
\[L(s,{J_{r,s,t;\delta}}_{/\BQ})=L^{(\infty)}(s-1/2,\phi_{r,s,t;\delta}),\]
where $\phi_{r,s,t;\delta}=\chi_\delta^{r+s}\phi_{r,s,t}$ is the Hecke character associated to the CM abelian variety ${J_{r,s,t;\delta}}_{/K}$  in the sense of \cite[Proposition 19.9]{Shimura-CM} (see \S 1.1) and $L^{(\infty)}$ means the $L$-function with the infinite factors removed. Therefore
\[L\left(s,{J_{r,s,t;\delta}}_{/\BQ}\right)=L^{(\infty)}(s-1/2,\phi_{r,s,t;\delta})=L^{(\infty)}(s-1/2,\pi_{r,s,t;\delta}).\]

Set the complete L-function as
\[\Lambda\left(s,{J_{r,s,t;\delta}}_{/\BQ}\right)=\N_{r,s,t;\delta}^{s/2}((2\pi)^{-s}\Gamma(s))^{(p-1)/2}L\left(s,{J_{r,s,t;\delta}}_{/\BQ}\right),\]
where $\N_{r,s,t;\delta}$ is the conductor of the curve ${C_{r,s,t;\delta}}_{/\BQ}$.
Then $\Lambda\left(s,{J_{r,s,t;\delta}}_{/\BQ}\right)$ satisfies Hecke's functional equation
\[\Lambda\left(s,{J_{r,s,t;\delta}}_{/\BQ}\right)=\epsilon(\pi_{r,s,t;\delta})\Lambda\left(2-s,{J_{r,s,t;\delta}}_{/\BQ}\right).\]
It follows immediately that
\[\epsilon(\pi_{r,s,t;\delta})=(-1)^{\ord_{s=1 }\Lambda\left(s,{J_{r,s,t;\delta}}_{/\BQ}\right)}.\]

By the theorem of Mordell-Weil, 
\[J_{r,s,t;\delta}(\BQ)\simeq \BZ^{R_{r,s,t;\delta}}\times J_{r,s,t;\delta}(\BQ)_\tor,\]
where $R_{r,s,t;\delta}\geq 0$ is an integer. Then the Birch and Swinnerton-Dyer conjecture predicts that the sign of the epsilon factor $\epsilon(\pi_{r,s,t;\delta})$ is compatible with the parity of the Mordell-Weil rank $R_{r,s,t;\delta}$, and we will discuss this issue in the final section.

\subsection{The $p$-torsions}
We determine the $p$-torsion points of the Jacobian varieties of Fermat curves over the cyclotomic field $K$ following the strategy presented in \cite{Greenberg80}.  We will first obtain the results for $\delta\in K^\times$, and then deduce the results for $\delta\in \BQ^\times$.

The field extension $K/\BQ$ is totally ramified at $p$.  Recall $\Pi=\omega-\ov{\omega}$, and  the ideal $\fP=(\Pi)$  is the unique prime ideal of $K$  above $p$.

Suppose $\delta\in K^\times/K^{\times p}$. For fixed $r,s,t$, we simply write $C_\delta$ for the curve
\[C_{r,s,t;\delta}:y^p=x^r(\delta-x)^s,\]
and $J_\delta$ for the Jacobian variety $J_{r,s,t;\delta}$. Consider the natural projection morphism of $K$-curves
\[\pr:C_{\delta}\ra \BP^1, \quad (x_0,y_0)\mapsto x_0,\]
which is associated to the natural embedding $K(x)\hookrightarrow K(x,y)$. Denote by
\[p_0=(0,1),\quad p_\delta=(\delta,1)\textrm{ and } p_\infty=(1,0)\]
the points on $\BP^1$ corresponding to the primes 
\[(x),\quad (\delta-x)\textrm{ and }(1/x)\]
of $K(x)$ respectively.  The morphism $\pr$ is totally ramified at these points $p_0,p_1$ and $p_\infty$. The corresponding fibers on $C_{\delta}$ are 
\[P_0=(0,0,1),\quad P_\delta=(\delta,0,1) \textrm{ and }P_\infty=(1,0,0)\]
respectively, which correspond to the primes 
\[(x,y),\quad (\delta-x,y) \textrm{ and }(1/x,y)\]
of $K(x,y)$ respectively. We also denote by $p_0,p_\delta,p_\infty, P_0,P_\delta,P_\infty,$  the corresponding prime ideals.

The rational divisor
\[D_0=P_0-P_\infty\]
is not principal and satisfies 
\[\div(x)=pD_0.\]

For any field $L\supset K$, denote by $\Div^0(C_\delta)(L)$ the group of $L$-rational divisor of degree zero on the curve $C_\delta$. Then the Mordell-Weil group $J_{\delta}(L)$ is identified with the quotient group of  $\Div^0(C_\delta)(L)$ by the principal divisors. Denote by $[D]$ the divisor class of $D$.
So the divisor class $[D_0]$ gives rise to a $p$-torsion point  in $J_{ \delta}(K)$. It is clear from the complex multiplication that $[D_0]\in J_{\delta}[\Pi]$. 
Since $J_{\delta}$ has complex multiplication by $\CO_K$, we have
\[J_{\delta}[\Pi^i]=\CO_K/(\Pi^i).\]
In particular, 
\[J_{ \delta}[\Pi]=\left\langle [D_0]\right\rangle.\]

The Galois group $G=\Gal(K(x,y)/K(x))$ is generated by the $K(x)$-automorphism 
\[\iota:K(x,y)\ra K(x,y),\quad \iota(y)=\omega y.\] 
Then $\iota$ gives the automorphism $\tau_2$ on $C_{\delta}$ and then $G$ acts on the group of divisors on $C_{\delta}$. Denote
\[\N=\sum_{i=0}^{p-1}\iota^i\textrm{ and }\mathrm{D}=1-\iota.\]
\begin{prop}\label{Hilbert90}
Let $L$ be a finite extension of $K$.  Then
\[\RH^1(G,\Div^0(C_\delta)(L))=0.\]
\end{prop}
\begin{proof}
Since $G$ is a cyclic group of order $p$, we have
\[\RH^1(G,\Div^0(C_\delta)(L))=\Div^0(C_\delta)(L)[N]/\mathrm{D}(\Div^0(C_\delta)(L)).\] 
In order to prove the proposition, we shall prove that if $D$ is an $L$-rational divisor with norm zero, i.e. $\N(D)=0$, then $D\in \mathrm{D}(\Div^0(C_\delta)(L))$. 

Let $D$ is an $L$-rational divisor with norm zero, i.e. $\N(D)=0$. Suppose
\[D=\sum_{s\in \BP^1}D_s\]
where $D_s$ is a divisor supported on the fiber $\pr^{-1}(s)$. Then each $D_s$ has norm zero and hence degree zero. Therefore, we may assume that $D=D_s$ is supported on a single fibre $\pr^{-1}(s)$ for some $s\in \BP^1(L)$. Note $G$ acts transitively on the fiber $\pr^{-1}(s)$. Suppose $S_0\in \pr^{-1}(s)$. Then
\[D=\sum_{S_i\in \pr^{-1}(s)} a_i S_i=\sum a_i (S_i-S_0)=\sum a_i(\iota^{r_i}-1)S_0=(1-\iota)D_1\]
where $\sum a_i=0$ and $S_i=\iota^{r_i}(S_0)$. Next we claim $D_1$ is $L$-rational. Since $D$ is $L$-rational, for any $\sigma\in G_L$, 
\[(1-\iota)(\sigma(D_1)-D_1)=0.\]
Since the only fixed points of $\iota$ on $C_{\delta}$ are $P_0$, $P_\delta$ and $P_\infty$, there exist integers $l,m,n\in \BZ$ with $l+m+n=0$ such that
\[\sigma(D_1)-D_1=lP_0+ mP_\delta+nP_\infty.\]
Suppose 
\[D_1=rP_0+sP_\delta+tP_\infty+D_2\]
where $D_2$ is not supported on $P_1,P_\delta,P_\infty$. Then
\[\sigma(D_1)-D_1=\sigma(D_2)-D_2=lP_0+ mP_\delta+nP_\infty.\]
We conclude that $l=m=n=0$. Thus $D_1$ is $L$-rational and $D\in \mathrm{D}(\Div^0(C_\delta)(L)).$
\end{proof}

\begin{prop}\label{torsion}
If $\delta\in K^{\times p}$, then $J_{r,s,t;\delta}(K)_\tor$ contains $J_{r,s,t;\delta}[\Pi^2]$. Otherwise, the $\Pi$-primary component of $J_{r,s,t;\delta}(K)_\tor$ is $J_{r,s,t;\delta}[\Pi]$. 
\end{prop}
\begin{proof}
It suffices to determine whether there is a divisor class $c_1$ such that 
\[(\iota-1)c_1=[D_0], \quad D_0=P_0-P_\infty.\] 
Such a divisor class $c_1$ exists if and only if there exists $z\in K(x,y)$ such that $x/\delta=\N(z)$.  By Proposition \ref{Hilbert90}, the existence of $c_1$ amounts to the existence of a divisor $D_0'\in [D_0]$ such that $\N(D_0')=0$.  If $\N(z)=x/\delta$ and take $D_0'=D_0-\div(z)$, then
\[\N(D_0')=\N(D_0)-\div(x)=0.\]
Conversely, if such a $D_0'$ exists, then there exists $z'\in K(x,y)$ with $\N(z')=a x$ with $a\in K^\times $. By considering the residue classes mod $P_\delta$, we see $a\delta=b^p$ for some $b\in K^\times$ and $z=z'/b$ has the required property.

Now by \cite[Lemma 1]{Greenberg80}, $x/\delta$ is a norm for the extension $K(x,y)/K(x)$ if and only if $x^r(\delta-x)^s$ is a norm for the extension $K(\sqrt[p]{x/\delta})/K(x)$. Note both $x/\delta$ and $1-x/\delta$ are norms for the extension $K(\sqrt[p]{x/\delta})/K(x)$. Hence
\[\delta^{-p}x^r(\delta-x)=\delta^{r+s-p}(x/\delta)^r(1-x/\delta)^s\]
is a norm if and only if $\delta$ is a norm, i.e. $\delta\in K^{\times p}$.
\end{proof}
\begin{remark}
If $\delta=1$ and $p\geq 5$, Greenberg \cite{Greenberg80} even proves $J_{r,s,t;1}[\Pi^3]\subset J_{r,s,t;1}(K)_\tor$.
\end{remark}
\begin{coro}
If $\delta$ is a nonzero rational number which is not in $ \BQ^{\times p}$, then $J_{r,s,t;\delta}(K)[\Pi^\infty]=J_{r,s,t;\delta}[\Pi]$.
\end{coro}
\begin{proof}
We show that $K^{\times p}\cap \BQ^\times=\BQ^{\times p}$. Consider the inflation-restriction sequence
\[0\ra\RH^1(\Gal(K/\BQ),\mu_p)\ra \RH^1(\BQ,\mu_p)\ra \RH^1(K,\mu_p).\]
Since $\Gal(K/\BQ)$ has cardinality $p-1$ prime to $p$, $\RH^1(\Gal(K/\BQ),\mu_p)=0$. Then via Hilbert satz 90, we have an injection
\[\BQ^\times/\BQ^{\times p}\hookrightarrow K^\times/K^{\times p}.\]
Then  we get $K^{\times p}\cap \BQ^\times=\BQ^{\times p}$, and the corollary follows by Proposition \ref{torsion}.

\end{proof}

\section{Root Numbers}
Suppose $\delta\in \BQ^\times/\BQ^{\times p}$ is an integer.  Fix the triple $(r,s,t)$ of positive integers satisfying $r+s+t=p$.  Simply denote $J_\delta=J_{r,s,t; \delta}$, and denote  $\phi=\phi_{r,s,t}$ and $\phi_\delta=\phi_{r,s,t;\delta}=\phi\chi_{\delta}^{r+s}$ the Hecke characters associated to $J$ and $J_{\delta}$. We shall compute the local root numbers of the Hecke character $\phi_\delta$. 
\subsection{Local root numbers}
We shall recall several formulae for the local root numbers. The basic references to this subject are \cite{Tate-thesis, JL70, Deligne73, Tate79}. Let $L$ be a nonarchimedean  local field. Fix a nontrivial additive character  $\psi$ of $L$, and denote by $dx$ the self-dual Haar measure on $L$ with respect to $\psi$. Let $\chi$ be a unitary character of $F^\times$. The local epsilon factor $\epsilon(\chi,\psi)=\epsilon(\chi,\psi,dx)$ associated to $\chi$ with respect to $\psi$ is given by \cite[(3.3.1)]{Deligne73}, and define
\[\epsilon(s,\chi,\psi)=\epsilon(\chi|\ |^s,\psi).\]
The local root number associated to $\chi$ with respect to $\psi$ is 
\[\frac{\epsilon(\chi,\psi)}{|\epsilon(\chi,\psi)|}=\epsilon(1/2,\chi,\psi).\]

If $L$ is archimedean, the formulae for the local epsilon factors can be found in \cite[(3.4.1)-(3.4.2)]{Deligne73} or \cite[p. 97]{JL70}.  Now suppose $L$ is nonarchimedean with ring of integers $\CO$, and a uniformizer $\varpi$.  Define the conductors: 
\begin{itemize}
\item[--] $n(\psi)=$ the largest integer $n$ such that $\psi|_{\varpi^{-n}\CO}=1$;
\item[--] $c(\chi)=$ the smallest positive integer $m$ such that $\chi|_{(1+\varpi^m\CO)}=1$ if $\chi$ is ramifed, and $0$ if $\chi$ is unramifed.
\end{itemize}
We also may refer to $\varpi^{-n}\CO$ resp. $\varpi^m\CO$ as the conductor of $\psi$ resp. $\chi$.

If $\chi$ is unramified, then 
\begin{equation}\label{unramified-root-number}
\epsilon(s,\chi,\psi)=\chi(\varpi^{n(\psi)})|\varpi|^{n(\psi)(s-\frac{1}{2})}.
\end{equation}
So if $n(\psi)=0$, then the root number 
\[\epsilon(1/2,\chi,\psi)=1.\]
If $\chi$ is ramified, then from \cite[(3.4.3.2)]{Deligne73} we have
\begin{equation}\label{ramified-root-number}
\epsilon(s,\chi,\psi)=\int_{\varpi^{-(n(\psi)+c(\chi))}\CO^\times}\chi^{-1}(x)|x|^{-s}\psi(x)dx.
\end{equation}

Suppose $M$ is a finite extension of $\BQ_\ell$, and let $L$ be a separable quadratic extension of $M$. This determines uniquely an element $\Delta$ in $L^\times$ modulo $M^\times$ such that
\begin{equation}\label{qd}
L=M[\Delta], \quad \Delta^2\in M^\times.
\end{equation}
Let $\psi_\ell:\BQ_p\ra \BC^\times$ be the additive character given by 
\begin{equation}\label{character-ell}
\psi_\ell(x)=e^{2\pi i \iota(x)}
\end{equation}
 where $\iota:\BQ_p/\BZ_p\hookrightarrow \BQ/\BZ$ is the natural embedding. Then define $\psi_L:L\ra\BC^\times$ by
\[\psi_L(x)=\psi_\ell(\Tr_{L/\BQ_\ell}(x)).\]
\begin{prop}[Fr\"ochlich-Queyrut {\cite[Theorem 3]{FQ73}}]\label{FQ}
Let $L$ be a separable quadratic extension of $M$, with $\Delta$ satisfying $(\ref{qd})$, and let $\chi$ be a character of $L^\times$ so that $\chi|_{M^\times}=1$. Then
\[\epsilon(1/2,\chi,\psi_L)=\chi(\Delta).\]
\end{prop}

\subsection{Kummer extensions}
The field extension $K/\BQ$ is totally ramified at $p$. 
Put 
\[\Pi=\omega-\ov{\omega} \textrm{ and }\pi=\Pi^2=\omega^2+\ov{\omega}^2-2.\]
Then the ideal $\fP=(\Pi)$ resp. $\fp=(\pi)$ is the unique prime ideal of $K$ resp. $F$ above $p$.

Note under the class field theory, the character $\chi_ \delta$ can be viewed as a character on the Galois group $\chi_\delta:G_K\ra \CO_K^\times$ which is characterized by $\chi_ \delta(\sigma)=({\sqrt[p]{ \delta}})^{\sigma-1}$.
\begin{lem}
\label{cond}
Let $\fD_0$ be the product of all distinct prime factors of $(\delta)$ in $K$ not dividing $p$. Then
\[\chi_{\delta}|_{\BA^\times}=1,\]
and the conductor $c(\chi_\delta)$ divides $\fP^{p+1} \fD_0$.
\end{lem}
\begin{proof}
Denote $L=K(\sqrt[p]{\delta})$.  We first note that the extension $L/K$ is unramified outside the places dividing $p\delta$. Under the Artin map, we have the isomorphism
\[\sigma:K^\times \N_{L/K}(\BA_L^\times)\backslash\BA_K^\times\simeq \Gal(L/K)\]
Let $c\in \Gal(K/F)$ be the non-trivial element. Then for any $\sigma\in \Gal(L/K)$, we have 
\[c\sigma_t c=\sigma_{t^c}=\sigma_t^{-1}.\] Let $v$ be a finite place of $F$.  For any $t\in K_v^\times$,
\[\sigma_{\N(t)}=\sigma_t\sigma_{t^c}\]
acts trivially on $\sqrt[p]{\delta}$, i.e. $\chi_{\delta}|_{\N(K_v^\times)}=1$.  If $v$ is split in $K$, then $\N(K_v^\times)=F_v^\times$. If $v$ is inert in $K$, then by class-field theory,  $\N(K_v^\times)$ is a subgroup of $F_v^\times$ of index $2$.  Since $\chi_\delta$ is a character of order $p$, in both cases,  $\chi_\delta$ is trivial when restricted on $F_v^\times$. By noting that $\chi_\delta$ has trivial local factors at infinite places, we have $\chi_\delta|_{\BA^\times}=1$.

For any $V\neq\Pi$ dividing $\delta$ and  $t\in (1+\fD_0\CO_{K,v})$,
\[\chi_\delta(t)=\left(\frac{t,\delta}{K_{V};p}\right)\equiv t^{- \frac{\ord_{V}(\delta)\left(q_{V}-1\right)}{p}}\equiv 1\mod \fD_0\CO_{K,V},\]
where $\left(\frac{\cdot,\cdot}{K_{V};p}\right)$ denotes the $p$-th Hilbert symbol over $K_V$.

For any $t=1+x\in K_\Pi^\times$ with $\ord_{\Pi}(x)>p$, it is easy to see that the binomial expansion
\[(1+x)^{\frac{1}{p}}=\sum_{i=0}^\infty \left(\frac{1/p}{i}\right)x^{i}\]
converges, and hence $t$ is a $p$-th power in $K_\Pi^\times$. Then
\[\chi_\delta(t)=\left(\frac{t,\delta}{K_{\Pi};p}\right)=\hilbert{\delta}{t}{K_\Pi}{p}^{-1}=1.\]

\end{proof}

\begin{lem}\label{ramification}
Let $\ell\neq p$ be a rational prime number. Let $V\mid \ell$ be a place of $K$. Then $V$ is fixed by the complex conjugation, or equivalently, $V$ is inert in the extension $K/F$,  if and only if $\ell$ has even order in $\BF_p^\times$. In particular, the places of $F$ above $\ell$ have the same inertia degree in the quadratic extension $K/F$.
\end{lem}\label{decom}
\begin{proof}
We identify the Galois group $\Gal(K/\BQ)$ with $(\BZ/p\BZ)^\times$ by the Artin reciprocity map, via which the complex conjugation is identified with $-1$.  Let $D_V$ be the decomposition group of $V$ with respect to the extension $K/\BQ$. Since $V$ is unramified, 
\[D_V\simeq \Gal(\BF_{q_V}/\BF_\ell).\]
Then $V$ is fixed by the complex conjugation if and only if $D_V$ contains $-1$. This is equivalent to that the inertia degree $[q_V:\ell]$ is even. It is well-known that the inertia degree $[q_V:\ell]$ is the order of $\ell$ in $\BF_p^\times$.

Since $K/\BQ$ is abelian, if any $D_V$ contains $-1$, so do the others, and the last assertion follows.
\end{proof}

 We write, in $\BQ_p^\times$, that
\[\delta=\epsilon p^a(1-p)^b\]
with $\epsilon\in \BF_p^\times\subset \BQ_p^\times$, $a\in \{0,1\cdots, p-1\}$, and $b\in \BZ_p$. Note such $a$ and $b$ are uniquely determined by
\[a=\ord_p(\delta),\quad b=\frac{\log(\delta)}{\log(1-p)},\]
where $\log$ denotes the $p$-adic Iwasawa logarithm (i.e. $\log(p)=0$). If $a=0$, we also have
\[\ord_p(b)=\ord_p\left(\frac{\delta^{p-1}-1}{p}\right).\]
Put
\begin{equation}\label{value-conductor}
u(\delta)=\min(\ord_p(a),\ord_p(b)+1).
\end{equation}
We note that $u(\delta)=\ord_p(b)+1$ if $a=0$, and $u(\delta)=0$ otherwise.
\begin{prop}\label{conductor}
The conductor of the Hecke character $\chi_\delta$ is given as
\[c(\chi_\delta)=\left\{
\begin{aligned}
&\fP^{p+1}\fD_0,& \quad \textrm{if  }u(\delta)=0;\\
&\fP^2\fD_0,&\quad  \textrm{if  }u(\delta)=1;\\
&\fD_0,&\quad  \textrm{if  }u(\delta)\geq 2;\\
\end{aligned}\right.\]
\end{prop}
\begin{proof}
This follows from Lemma \ref{cond} and \cite[Theorem 6.1]{CM88} or \cite[Theorem 8]{Sharifi01}.
\end{proof}

Put 
\[\Pi'=1-\omega.\]
We have
\[\Pi=-2\Pi'+\Pi'\ov{\Pi'}, \quad \omega=1-\Pi',\]
and note 
\begin{equation}\label{FQ-condition}
K=F[\Pi],\quad \Pi^2\in F^\times.
\end{equation}
For any place $V$ of $K$, denote by $\chi_{\delta,V}$ the $V$-adic component of $\chi$. By Proposition \ref{conductor}, if $u(\delta)\geq 2$, then the local character $\chi_{\delta,\Pi}$ is unramified.
\begin{prop}\label{value-conductor-chi}
\begin{itemize}
\item[1.] Suppose $u(\delta)=1$. The local character $\chi_{\delta,\Pi}$ has conductor $(\Pi^2)$, and 
\[\chi_{\delta,\Pi}(1+\Pi)=\omega^{2\left(\frac{ \delta^{p-1}-1}{p}\right)}.\]

\item[2.] Suppose $u(\delta)=0$. The local character $\chi_{\delta,\Pi}$ has conductor $(\Pi^{p+1})$, and 
\[\chi_{\delta,\Pi}(1+\Pi^{p})=\omega^{2\ord_p(\delta)}.\]
\end{itemize}

\end{prop}
\begin{proof}
First suppose $u(\delta)=1$. In this case, $a=0$ and $\chi_{\delta,\Pi}$ has conductor $(\Pi^2)$.
Since $\epsilon$ is a $p$-th power in $K_\Pi^\times$, we have
\[\chi_{\delta,\Pi}(1+\Pi)=\hilbert{1+\Pi}{\delta}{K_\Pi}{p}=\hilbert{1+\Pi}{(1-p)^b}{K_\Pi}{p}.\]
By an explicit reciprocity formula of Artin-Hasse-Iwaswa \cite{Iw68}, we have
\[\hilbert{1+\Pi}{(1-p)^b}{K_\Pi}{p}=\omega^{[1+\Pi, (1-p)^b]}\]
where
\[[1+\Pi, (1-p)^b]=-\frac{1}{p}\Tr_{K_\Pi/\BQ_p}\left(\frac{\omega}{1+\Pi}\frac{d (1+\Pi)}{d\Pi'}\log (1-p)^b\right).\]
Note
\[\log(1-p)^b=-bp\left(1+\frac{p}{2}+\frac{p^2}{3}+\cdots\right).\]
\[\]
And
\begin{eqnarray*}
\frac{\omega}{1+\Pi}\frac{d(1+\Pi)}{d\Pi'}&=& \frac{1-\Pi'}{1-2\Pi'+\Pi'\ov{\Pi'}}\cdot \frac{d(1-2\Pi'+\Pi'\ov{\Pi'})}{d\Pi'}\\
&\equiv& -2 \mod \Pi'.
\end{eqnarray*}
Hence
\[-\frac{1}{p} \Tr_{K_\Pi/\BQ_p}\left( \frac{\omega}{1+\Pi} \frac{d(1+\Pi)}{d\Pi'}\log(1-p)^b\right)\equiv -\Tr_{K_\Pi/\BQ_p}(2b) \equiv 2b\mod p.\]
On the other hand, 
\[\delta^{p-1}=(\epsilon(1-p)^b)^{p-1}\equiv 1+bp \mod p^2.\]
So
\[-\frac{1}{p}\Tr_{K_\Pi/\BQ_p}\left( \frac{\omega}{1+\Pi} \frac{d(1+\Pi)}{d\Pi'}\log(1-p)^b\right)\equiv 2\left(\frac{\delta^{p-1}-1}{p}\right)\mod p.\]
The first assertion follows.

Next suppose $u(\delta)=0$.  In this case, $a\neq 0$, and the local character $\chi_{\delta,\Pi}$ has conductor $(\Pi^{p+1})$. We have
\[\chi_{\delta,\Pi}(1+\Pi^p)=\hilbert{1+\Pi^p}{p^a(1-p)^b}{K_\Pi}{p}.\]
By Proposition \ref{conductor}, the local character $\chi_{(1-p)^b,\Pi}$ has conductor dividing $(\Pi^2)$. Hence
\[\chi_{\delta,\Pi}(1+\Pi^p)=\hilbert{1+\Pi^p}{p^a}{K_\Pi}{p}.\]
Since
\[1+\Pi^p=1+(-2\Pi'+\Pi\ov{\Pi'})^p\equiv 1+(-2)^p\Pi'^p\equiv(1+\Pi'^p)^{-2}\mod \Pi^{p+1},\]
\[\chi_{\delta,\Pi}(1+\Pi^p)=\hilbert{1+\Pi'^p}{p}{K_\Pi}{p}^{-2a}.\]
It suffices to compute 
\[\hilbert{1+\Pi'^p}{p}{K_\Pi}{p},\]
and we will apply Coleman's explicit reciprocity law \cite[Theorem 1]{Coleman81} and \cite[p. 89]{CM88} to do this calculation. Let $g_1(T)\in \BZ_p[[T]]$ be as in p. 90 of \cite{CM88}:
\[g_1(T)=\prod_{a=1}^{p-1}[a](T),\]
where $[a](T)=1-(1-T)^a$ be the  $a$-th power homomorphism of the multiplicative formal group $\wh{\BG}_\m$ with group law $\wh{\BG}_\m(X,Y)=X+Y-XY$. Then
\[g_1(\Pi')=p.\]
Define
\[\Lambda:1+T\BZ_p[[T]]\lra T\BZ_p[[T]],\quad \Lambda(f)=\log(f)-\frac{1}{p}\log(f([p])).\]
By Coleman's explicit reciprocity law
\[\hilbert{1+\Pi'^p}{p}{K_\Pi}{p}=\omega^{-n},\]
where
\[n=\int_1 \Lambda(1+T^p)\frac{Dg_1}{g_1},\]
and we note the Hilbert symbol is the inverse of that in Coleman's paper \cite[p. 43]{CM88}.
It is straight-forward to verify
\[\Lambda(1+T^p)\equiv (1-p^{p-1})T^p\mod T^{p+1}\BZ_p[[T]].\]
Since $\frac{Dg_1}{g_1}$ has at most a single pole at $T=0$, we obtain

\[n=\frac{1}{p}\Tr_{K_\Pi/\BQ_p}\left( \left.\left(\Lambda(1+T^p)\frac{Dg_1}{g_1}\right)\right|_{T=\Pi'}\right).\]

Since $\chi_{p,\Pi}$ has conductor $p+1$, by \cite[Theorem 6.3]{CM88},  we have
\[n\equiv \frac{1}{p}\Tr_{K_\Pi/\BQ_p}\left( \left.\left((1-p^{p-1})T^p\frac{Dg_1}{g_1}\right)\right|_{T=\Pi'}\right)\equiv \frac{1}{p} \Tr_{K_\Pi/\BQ_p}\left( \left.\left(T^p\frac{Dg_1}{g_1}\right)\right|_{T=\Pi'}\right)\mod p.\]
It follows from \cite[Lemma 6.5]{CM88} that (There is a sign error in this lemma: indeed, the expression for $D[a]/[a]$ in \cite[p. 91]{CM88} has a sign error.)
\[n\equiv\frac{1}{p}\Tr_{K_\Pi/\BQ_p}(\Pi'^{p-1}) \mod p.\]
On the other hand,
\begin{eqnarray*}
\Tr_{K_\Pi/\BQ_p}(\Pi'^{p-1}) &=&\sum_{j=1}^{p-1}(1-\omega^j)^{p-1}\\
&=&p-1+\sum_{i=1}^{p-1}\begin{pmatrix}p-1\\i\end{pmatrix}(-1)^i(\omega^i+\omega^{2i}+\cdots+\omega^{(p-1)i})\\
&=&p-1+\sum_{i=1}^{p-1}\begin{pmatrix}p-1\\i\end{pmatrix}(-1)^{i+1}\\
&=&p-\sum_{i=0}^{p-1}\begin{pmatrix}p-1\\i\end{pmatrix}(-1)^{i}=p.
\end{eqnarray*}
So 
\[n\equiv 1\mod p,\]
and hence
\[\chi_{\delta,\Pi}(1+\Pi^{p})=\omega^{2a}.\]
\end{proof}

For any prime $\ell$, denote 
\[\chi_{\delta,\ell}=\prod_{V\mid \ell} \chi_{\delta,V}\]
where $V$ runs through the places of $K$ above $\ell$.

\begin{prop}\label{value-chi}
For any prime $\ell$, we have $\chi_{\delta,\ell}(\Pi)=1$.

\end{prop}

\begin{proof}
If $\ell=\infty$, this is clear because the infinite factor $\chi_{\delta,\infty}=1$. Next we assume $\ell$ is finite. By our choice of $\Pi$, we have $\Pi^2\in F$. By Lemma \ref{cond}, 
\[\chi_{\delta,\ell}(\Pi^2)=1, \textrm{ i.e. } \chi_{\delta,\ell}(\Pi)=\pm 1.\] On the other hand, by definition, we know $\chi_\delta(\Pi)$ is a $p$-th root of unity, and we conclude 
\[\chi_{\delta,\ell}(\Pi)=1.\]

\end{proof}

\subsection{Local root numbers outside $p$}
For any finite place $v$ resp. $V$ of $F$ resp. $K$ outside $p$, denote by $\varpi_v$ resp. $\varPi_V$ the corresponding uniformizers of $F_v$ resp. $K_V$.

We fix an additive character $\psi_{K_V}$ of $K_V$ for each place $V$ of $K$  as follows: If $V$ is archimedean,  then 
\[\psi_{K_V}(x)=e^{-2\pi i \Tr_{K_V/\BR}(x)};\] 
If $V$ is non-archimedean and dividing $\ell$, then 
\[\psi_{K_V}(x)=\psi_\ell\left(\Tr_{K_V/\BQ_\ell}(x)\right),\] 
where $\psi_\ell$ is as in (\ref{character-ell}). Then the product $\psi=\prod_V\psi_{K_V}$ is an additive character on $K\backslash \BA_K$.  The L-function $L(s,\phi_\delta)$ satisfies a functional equation
\[L(s,\phi_\delta)=\epsilon(s,\phi_\delta)L(1-s,\phi_\delta^{-1}).\]
The epsilon-factor admits a local product
\[\epsilon(s,\phi_\delta)=\prod_V\epsilon(s,\phi_{\delta,V},\psi_{K_V}),\]
and the epsilon-factor $\epsilon(s,\phi_\delta)$ is independent of the choices of $\psi_{K_V}$.

Let $V$ be an infinite place of $K$. Then
\[K_V\simeq \BC, \quad \phi_{\delta,V}(z)=z^{-1}|z|^{1/2}_\BC.\]
and, from \cite[(3.3.3), (3.4.2)]{Deligne73} or \cite[p. 97]{JL70}, we have
\begin{equation}\label{rn-inf}
\epsilon(1/2,\phi_{\delta,V},\psi_{K_V})=i^{-1}.
\end{equation}
For a place $V\nmid p\delta $ of $K$, $\phi_{\delta,V}$ is unramified, and  we have by (\ref{unramified-root-number}) that
\begin{equation}\label{rn-unr}
\epsilon(1/2,\phi_{\delta,V},\psi_{K_V})=1.
\end{equation}

Let $V$ be a place where $\phi_{\delta,V}$ is ramified. Simply denote $c_V=c(\phi_{\delta,V})$ and $n_V=n(\psi_{K_V})$. It follows from (\ref{ramified-root-number}) that 
\begin{equation}\label{epsilon1}
\epsilon(1/2,\phi_{\delta,V},\psi_{K_V})=\phi_{\delta,V}(\varPi_V)^{n_V+c_V}\N_V ^{\frac{n_V-c_V}{2}}\int_{\CO_{K,V}}dx \sum_{a\in \CO_{K,V}^\times/(1+\varPi_V^{c_V}\CO_{K,V})} \phi_{\delta,V}^{-1}(a)\psi_{K_V}\left(\frac{a}{\varPi_V^{n_V+c_V}}\right).
\end{equation}
Since $\phi$ is unramified outside $p$, it  follows immediately from the above formula that
\begin{lem}\label{rn1}
For any place $V\nmid p$, we have
\[\epsilon(1/2,\phi_{\delta,V},\psi_{K_V})=\epsilon(1/2,\chi_{\delta,V},\psi_{K_V}) \phi_{V}(\varPi_V)^{n_V+c_V}.\]
\end{lem}
\begin{proof}
This lemma follows from $(\ref{epsilon1})$ by noting that $\phi_\delta=\phi\chi_\delta^{r+s}$ and that $\phi$ is unramified at the place $V$.
\end{proof}
In order to compute the local root numbers for $\phi_{\delta,V}$ for $V\nmid p$ and $V\mid \delta$, it suffices to compute $\phi_V(\varPi_V)$ and $\epsilon(1/2,\chi_{\delta,V},\psi_{K_V})$ for these places.
Since our choice of $\Pi$ satisfies (\ref{qd}), by Proposition \ref{FQ}, for any finite place $V$ of $K$,
\begin{equation}\label{epsilon-chi}
\epsilon(1/2,\chi_{\delta,V},\psi_{K_V})=\chi_{\delta,V}(\Pi)=1,
\end{equation}
which is given by Proposition \ref{value-chi}. So it remains to compute the values of $\phi_V(\varPi_V)$ for $V\neq \Pi$.
\begin{prop}\label{value-phi1}
Let $v\neq \pi $ be a finite place of $F$. If $v$ is inert in $K$, let $V$ be the unique place above $v$.
Then
\[\phi_V(\varPi_V)=-1.\]
If $v$ is split in $K$, let $V$ and $\ov{V}$ be two distinct places above $v$. Then
\[ \phi_V(\varPi_V)\phi_{\ov{V}}(\varPi_{\ov{V}})=+1.\]
\end{prop}
\begin{proof}
First suppose $v$ is inert in $K$. The prime ideal $V$ of $K$ is fixed by the complex conjugation. By the Stickelberger relation (\ref{Stic}) and the fact that $\Phi_{r,s,t}$ is a CM-type, i.e. a set of coset representatives of $(\BZ/p\BZ)^\times$ modulo $\{\pm 1\}$, $j_{r,s,t}(V)$  generates the same ideal as $q_v=q_V^{1/2}$. Hence 
\[j_{r,s,t}(V)=u q_v\]
for some unit $u\in \CO_K^\times$. Since $u$ has absolute value $1$ for any complex embedding, it must be a root of unity. Note $q_v\equiv -1\mod p$. Then the congruence
\[ j_{r,s,t}(V)\equiv 1\mod \Pi^2\]
forces $u=-1$. Then 
\[\phi_V(\varPi_V)=j_{r,s,t}(V)\cdot q_V^{-1/2}=-1.\]

Next, we suppose $v$ is split in $K$. Then by definition we have
\[\phi_V(\varPi_V)\equiv -\sum_{\begin{subarray}{c}x\in \CO_{K}/V\\x \neq 0,1\end{subarray}} \left(x^r(1-x)^s\right)^{\frac{q_V-1}{p}}\mod V,\]
and 
\[\phi_{\ov{V}}(\varPi_{\ov{V}})\equiv -\sum_{\begin{subarray}{c}x\in \CO_{K}/\ov{V}\\x \neq 0,1\end{subarray}} \left(x^r(1-x)^s\right)^{\frac{q_{\ov{V}}-1}{p}}\mod \ov{V}.\]
 Since $v$ is split in $K$, we have natural identities
\[\CO/v=\CO_K/V=\CO_K/\ov{V},\]
and we may take the representatives in $\CO$ for each $x$.
Taking complex conjugation, we have then
\[\ov{\phi_{\ov{V}}(\varPi_{\ov{V}})}\equiv -\sum_{\begin{subarray}{c}x\in \CO_{K}/\ov{V}\\x \neq 0,1\end{subarray}} \left(x^r(1-x)^s\right)^{\frac{q_v-1}{p}}\mod V.\]
So we have 
\[\ov{\phi_{\ov{V}}(\varPi_{\ov{V}})}=\phi_V(\varPi_V),\]
or equivalently
\[ \phi_V(\varPi_V)\phi_{\ov{V}}(\varPi_{\ov{V}})=+1\]
\end{proof}
\begin{coro}\label{value-phi2}
Denote
\[\phi_\ell=\prod_{V\mid \ell}\phi_V(\varPi_V)\]
for each prime $\ell\neq p$. Then
\[
\phi_\ell=\left(\frac{\ell}{p}\right).
\]

\end{coro}
\begin{proof}
Let $f$ be the order of $\ell$ in $\BF_p^\times$. Then $f$ is the inertia degree of the places of $K$ above $\ell$. If $\left(\frac{\ell}{p}\right)=+1$, then $(p-1)/f$ is even, i.e. there are even places of $K$ above $\ell$. Then it follows from Lemma \ref{ramification} and Proposition \ref{value-phi1} that  if all the places of $F$ above $\ell$ are split in $K$, 
\[\phi_\ell=(+1)^{(p-1)/2f}=+1,\]
and if all the places of $F$ above $\ell$ are inert in $K$,
\[\phi_\ell=(-1)^{(p-1)/f}=+1.\]

If $\left(\frac{\ell}{p}\right)=-1$, then $(p-1)/f$ is odd. By Lemma \ref{ramification}, the places of $K$ above $\ell$ must be inert in the quadratic extension $K/F$. Then in this case we have 
\[\phi_\ell=(-1)^{(p-1)/f}=-1.\]
\end{proof}


Recall we simply denote $c_V=c(\phi_{\delta,V})$ and $n_V=n(\psi_{K_V})$ for any place $V$ of $K$.
\begin{lem}\label{conductor1}
\begin{itemize}
\item[1.] The conductors $n_\Pi=1$ and $n_V=0$ for all $V\neq \Pi$.
\item[2.] The conductors $c_V=0$ if $V\nmid p\delta$, $c_V=1$ if $V\mid \delta$ but $V\neq \Pi$, and \[
c_\Pi=\left\{\begin{aligned}
&\leq 2,&\quad \textrm{ if } p\nmid \delta;\\
&p+1,&\quad \textrm{ if } p\mid \delta.
\end{aligned}
\right.
\]
\end{itemize}
\end{lem}
\begin{proof}
Recall the relative different $\partial_{K/F}$ is defined to be the integral ideal of $K$ such that
\[\partial_{K/F}^{-1}=\{x\in K|\,\tr_{K/F}(xy)\in \CO_F\textrm{ for all }y\in \CO_F\}.\]
It follows from the discriminant formula
\[d_{K}=d_F^2\cdot \N_{K/\BQ}(\partial_{K/F})\]
that $\N_{K/\BQ}(\partial_{K/F})=(p)$, and hence $\partial_{K/F}=(\Pi)$. By the definition of $\psi_{K_V}$, $n_V$ is the exponent of $\partial_{K/F}$ at $V$. Then the values for $n_V$ follows. The second assertion follows from Proposition \ref{conductor} and the conductors of $\phi$.
\end{proof}


Let $\ell$ be a prime number, including the infinity place $\infty$. Denote 
\[\epsilon_\ell(\phi_{r,s,t;\delta})=\prod_{V\mid \ell} \epsilon(1/2,\phi_{r,s,t;\delta,V},\psi_{K_V)}\]
where $V$ runs through the places of $K$ above $\ell$.

\begin{prop}\label{lrn1}
Let $\ell\neq p$ be a prime number.
\begin{itemize}
\item[1.] If $\ell=\infty$ is archimedean, then
\[\epsilon_\infty(\phi_{r,s,t;\delta})=i^{-\frac{p-1}{2}};\]
If $\ell$ is a finite prime with $\ell\nmid p\delta$, then
\[\epsilon_\ell(\phi_{r,s,t;\delta})=1;\]

\item[2.] If $\ell\neq p$ and $\ell\mid \delta$, then
\[\epsilon_\ell(\phi_{r,s,t;\delta})=\left(\frac{\ell}{p}\right).\]
\end{itemize}
\end{prop}

\begin{proof}
The first assertion are just (\ref{rn-inf}) and (\ref{rn-unr}). By Lemma \ref{rn1} and Lemma \ref{conductor1}, 
\[\epsilon(1/2,\phi_{\delta,V},\psi_{K_V})=\epsilon(1/2, \chi_{\delta,V},\psi_{K_V})\phi_V(\varPi_V)\]
for any finite place $V\neq \Pi$ and $V\mid \delta$. Then the remaining follows from (\ref{epsilon-chi}) and Corollary \ref{value-phi2}.
\end{proof}

\subsection{Local root number at $p$}
Let $\eta$ be the quadratic character of $F^\times\backslash\BA^\times$ determined by the quadratic extension $K/F$ through class-field theory. Then $\eta$ is unramified outside $p$. For any finite place  $v$ of $F$ not dividing $p$, we have
\begin{equation}\label{epsilon-eta1}
\epsilon(1/2,\eta_v)=+1.
\end{equation}
For each archimedean place $\tau:F\hookrightarrow \BR$, we have
\[F_\tau\simeq \BR, \textrm{ and }\eta_\tau(x)=\frac{|x|}{x},\, x\neq 0.\]
Then we have by \cite[p. 97]{JL70} that
\begin{equation}\label{epsilon-eta2}
\epsilon(1/2,\eta_\tau,\psi_{\BR})=i^{-1}.
\end{equation}
Finally, for the place $\pi$, we have
\begin{equation}\label{epsilon-eta}
\epsilon(1/2,\eta_\pi,\psi_{F_\pi})=i^{\frac{p-1}{2}}.
\end{equation}

The character $\phi$ is equivariant with respect to the complex conjugation. In particular, we have
\[\phi|_{\BA^\times}=\eta.\]
Since 
\[\chi_\delta|_{\BA^\times}=1,\]
we have
\[\phi_{\delta}=\phi\chi_\delta^{r+s}|_{\BA^\times}=\eta.\]

Recall for any $\delta\in \BQ^\times \cap \BZ_p$, we have the following decomposition in $\BQ_p^\times$:
\[\delta=\epsilon p^a(1-p)^b\]
with $\epsilon\in \BF_p^\times\subset \BQ_p^\times$, $a\in \{0,1\cdots, p-1\}$, and $b\in \BZ_p$.  We have defined in (\ref{value-conductor}) that
\[u(\delta)=\min(\ord_p(a),\ord_p(b)+1).\]
\begin{prop}\label{lrn2}
\begin{itemize}
\item[1.] If 
\[u(r^rs^s(t-p)^{t}\delta^{r+s})\geq2,\]
 then  $c(\phi_{\delta,\Pi})=1$, and 
 \[\epsilon(1/2,\phi_{\delta,\Pi},\psi_{K_\Pi})=\left(\frac{2}{p}\right) i^{\frac{p-1}{2}}.\]
 \item[2.] If 
 \[u(r^rs^s(t-p)^{t}\delta^{r+s})=1,\]
 then $c(\phi_{\delta,\Pi})=2$, and 
 \[\epsilon(1/2,\phi_{\delta,\Pi},\psi_{K_\Pi})=-\left(\frac{-rstd}{p}\right)i^{\frac{p-1}{2}},\]
where
\[d\equiv \frac{\left(r^rs^s(t-p)^t\delta^{r+s}\right)^{p-1}-1}{p}\mod p.\]
\item[3.]  If 
\[u(r^rs^s(t-p)^{t}\delta^{r+s})=u(\delta^{r+s})=0,\]
 then  $c(\phi_{\delta,\Pi})=p+1$, and 
 \[\epsilon(1/2,\phi_{\delta,\Pi},\psi_{K_\Pi})=-\left(\frac{rstd}{p}\right)i^{\frac{p-1}{2}},\]
 where 
 \[d\equiv (r+s)\ord_p(\delta)\mod p.\]
\end{itemize}
\end{prop}
\begin{proof}
First we clarify the conditions for the conductors of $\phi_{\delta,\Pi}$. By \cite[Theorem 5.3]{CM88}, for any $x\in \CO_{K,\Pi}^\times$, 
\begin{equation}\label{x}
\phi_{\delta,\Pi}(x)=\hilbert{x}{r^rs^s(t-p)^{t}\delta^{r+s}}{K_\Pi}{p}\hilbert{x}{\Pi}{K_\Pi}{2}.
\end{equation}
Then the statements on the conductors $c(\phi_{\delta,\Pi})$ follows from Proposition \ref{conductor} by noting that the character
\[x\mapsto \hilbert{x}{\Pi}{K_\Pi}{2}\]
has conductor $1$.

If $c(\phi_{\delta, \Pi})=1$, then by \cite[Proposition 2]{Rohr92}
\[\frac{\epsilon(1/2, \phi_{\delta,\Pi},\psi_{K_\Pi})}{\epsilon(1/2,\eta_\pi,\psi_{F_\pi})}=\left(\frac{2}{p}\right),\]
and the first assertion follows by (\ref{epsilon-eta}).

Again by \cite[Proposition 2]{Rohr92}, if $c(\phi_{\delta ,\Pi})=c>1$, then
\begin{equation}\label{roh2}
\frac{\epsilon(1/2, \phi_{\delta,\Pi},\psi_{K_\Pi})}{\epsilon(1/2,\eta_\pi,\psi_{F_\pi})}=\left(\frac{-2l}{p}\right)\phi_{\delta,\Pi}(\Pi)^{c-1}i^u,
\end{equation}
where $u$ is given by
\[u=\left\{
\begin{aligned}
0,\quad \textrm{ if } p\equiv 1\mod 4;\\
1,\quad \textrm{ if } p\equiv 3\mod 4,
\end{aligned}
\right.\]
and  $l$ is given by
\[\phi_{\delta,\Pi}(1+\Pi^{c-1})=\omega^l.\]

We compute the values of $\phi_{\delta,\Pi} (\Pi)=\phi_{\Pi}(\Pi)\chi_{\delta,\Pi}^{r+s}(\Pi)$. By Proposition \ref{value-chi}, it suffices to compute $\phi_{\Pi}(\Pi)$. Since $\phi$ is a Hecke character unramified outside $\Pi$ and infinity, we have
\[\phi_\Pi(\Pi)=\prod_{V\mid \infty} \phi_V(\Pi)^{-1}=\prod_{h\in \Phi^{-1}} \frac{\Pi^{\sigma_h}}{|\Pi^{\sigma_h}|_\BC^{1/2}}=\prod_{h\in \Phi^{-1}}i\cdot \sgn(\sin(2\pi h/p))=i^{\frac{p-1}{2}}(-1)^{\frac{p-1}{2}-k},\]
where $k$ denotes the number of $h\in \Phi^{-1}$ satisfying $1\leq h\leq (p-1)/2$.  By \cite[Proposition 4]{Rohr92},
\[(-1)^{k}=-\left(\frac{-2rst}{p}\right).\]

Next we determine the value of $l$. Suppose $c(\phi_{\delta,\Pi})=2$, or equivalently, $u(r^rs^s(t-p)^{t}\delta^{r+s})=1$.  By (\ref{x}),
\[\phi_{\delta, \Pi}(1+\Pi)=\hilbert{1+\Pi}{r^rs^s(t-p)^t\delta^{r+s}}{K_\Pi}{p}.\]
Then by Proposition \ref{value-conductor-chi}, 
\[\phi_{\delta,\Pi}(1+\Pi)=\omega^{2\left(\frac{\left(r^rs^s(t-p)^{t}\delta^{r+s}\right)^{p-1}-1}{p}\right)}.\]
In this case, we have
\[l\equiv {2\frac{\left(r^rs^s(t-p)^{t}\delta^{r+s}\right)^{p-1}-1}{p}}\mod p.\]
Then the formula in the second assertion follows by noting
\[i^{\frac{p-1}{2}+u}=\left(\frac{2}{p}\right) \textrm{ and }\epsilon(1/2,\eta_\pi,\psi_{F_\pi})=i^{\frac{p-1}{2}}.\]

Finally, suppose $c(\phi_{\delta,\Pi})=p+1\geq 4$. Since $\phi_\Pi$ has conductor $1$ or $2$,
\[\phi_{\delta,\Pi}(1+\Pi^{c-1})=\chi_{\delta,\Pi}^{r+s}(1+\Pi^{p}).\]
Then by Proposition \ref{value-conductor-chi},
\[\phi_{\delta,\Pi}(1+\Pi^{c-1})=\omega^{2(r+s)\ord_p(\delta)},\]
and in this case, 
\[l\equiv 2(r+s)\ord_p(\delta)\mod p.\]
Then the formula in the third assertion follows by noting 
\[i^{\frac{p(p-1)}{2}+u}=\left(\frac{-2}{p}\right) \textrm{ and }\epsilon(1/2,\eta_\pi,\psi_{F_\pi})=i^{\frac{p-1}{2}}.\]

\end{proof}

\subsection{Global root numbers}
Fix the triple $(r,s,t)$ of integers with $r,s,t>0$ and $r+s+t=p$. For any prime $\ell$ of $\BQ$, denote
\[\epsilon_\ell(\phi_{r,s,t;\delta})=\prod_{V\mid \ell } \epsilon(1/2,\phi_{r,s,t;\delta, V},\psi_{K_V}),\]
\[\epsilon_\ell(\pi_{r,s,t;\delta})=\prod_{v\mid \ell } \epsilon(1/2,\pi_{r,s,t;\delta, v},\psi_{F_v}),\]
where $V$ resp. $v$ runs through all places of $K$ resp. $F$ above $\ell$. Then by \cite[Lemma 1.2, Theorem 4.7]{JL70}, we have
\[\epsilon_\ell(\pi_{r,s,t;\delta})=\epsilon_\ell(\phi_{r,s,t;\delta})\epsilon_\ell(\eta).\]
Let $d\in (\BZ/p\BZ)^\times$ be such that
\[d\equiv\left\{\begin{aligned} 
\frac{\left(r^rs^s(t-p)^t\delta^{r+s}\right)^{p-1}-1}{p}\mod p,\quad &
\textrm{if $u(r^rs^s(t-p)^t\delta^{r+s})=1$}; \\
(r+s)\ord_p(\delta)\quad \quad\mod p,\quad &
\textrm{if $u(r^rs^s(t-p)^t\delta^{r+s})=0$.}\end{aligned}\right.\]
\begin{thm}\label{root-number}
The global root number of the automorphic representation $\pi_{r,s,t;\delta}$ is given as follows:
\[\epsilon(\pi_{r,s,t;\delta})=\prod_{\ell}\epsilon_\ell(\pi_{r,s,t;\delta}),\]
where 
\[\epsilon_\ell(\pi_{r,s,t;\delta})=\left\{\begin{aligned}
\left(\frac{-1}{p}\right),\quad &\textrm{if $\ell=\infty$};\\
\left(\frac{\ell}{p}\right),\quad & \textrm{ if $\ell\neq p$ and $\ell\mid \delta$;}\\
+1,\quad & \textrm{if $\ell\nmid p\delta$;}
\end{aligned}\right.\]
and 
\[\epsilon_p(\pi_{r,s,t;\delta})=\left\{
\begin{aligned}
-\left(\frac{-rstd}{p}\right),\quad &\textrm{ if } u(r^rs^s(t-p)^t\delta^{r+s})= 0;\\
-\left(\frac{rstd}{p}\right),\quad &\textrm{ if } u(r^rs^s(t-p)^t\delta^{r+s})= 1;\\
\left(\frac{-2}{p}\right),\quad &\textrm{ if } u(r^rs^s(t-p)^t\delta^{r+s})\geq 2.
\end{aligned}
\right.\]
\end{thm}
\begin{proof}
This is the combination of Propositions \ref{lrn1} and \ref{lrn2} and the formulae (\ref{epsilon-eta1}-\ref{epsilon-eta})
\end{proof}
If we take $\delta=1$, then this is the main result of \cite{Rohr92} for the case $n=1$.  Next suppose $p=3$, $r=s=t=1$ and $\delta \in \BQ^\times/\BQ^{\times 3}$. The Fermat curve
\[F_\delta:X^3+Y^3=\delta\] and the quotient curve $C_{1,1,1;\delta}$ are both elliptic curves  over $\BQ$ and $\varphi_{1,1,1}: F_\delta\ra C_{1,1,1;\delta}$ is an isogeny over $\BQ$. The elliptic curve $F_\delta$ has Weierstrass equation 
\[y^2=x^3-2^43^3\delta^2.\]
\begin{coro}
Let $\delta \in \BQ^\times/\BQ^{\times 3}$ be an integer. The global root number of the elliptic curve  $F_\delta$ is given by
\[\epsilon=\prod_{\ell\mid \infty 3 \delta}\epsilon_\ell\]
where 
\[\epsilon_{\infty}=-1;\]
and if $\ell\neq 3$ and $\ell\mid \delta$, then
\[\epsilon_\ell=\left(\frac{\ell}{3}\right);\]
and 
\[\epsilon_3=\left\{
\begin{aligned}
-1,\quad &\textrm{ if } \delta\equiv \pm 1\mod 9, \textrm{ or } \ord_3(\delta)=1;\\
+1,\quad &\textrm{ otherwise}.
\end{aligned}
\right.\]
\end{coro}
\begin{proof}
It suffices to translate the expressions of $\epsilon_3$ from  Theorem \ref{root-number}, since the expressions for $\epsilon_\ell$ are clear for $\ell \neq 3$. In this case, we have $r=s=t=1$. First suppose $u(r^rs^s(t-p)^t\delta^{r+s})= 0$ i.e. $a=\ord_3(\delta)>0$. Then 
\[\epsilon_3=-\left(\frac{-2a}{3}\right),\]
and $\epsilon_3=-1$ if and only if $a=1$. Next suppose $u(r^rs^s(t-p)^t\delta^{r+s})= 1$. Then
\[\epsilon_3=-\left(\frac{d}{3}\right),\]
and it is straight-forward to verify that $\epsilon_3=-1$ if and only if $\delta\equiv \pm 1\mod 9$. Finally, if $u(r^rs^s(t-p)^t\delta^{r+s})= 2$, then 
\[\epsilon_3=\left(\frac{-2}{3}\right)=+1.\]
Combining all these, the lemma follows.
\end{proof}

This is compatible with the results of  Birch-Stephens \cite{BS66,Stephens68} and Liverance \cite{Liver95} which gives the formulae for root numbers of the elliptic curves with Weierstrass equation 
\[y^2=x^3+D\]
for $D\in \BQ^\times/\BQ^{\times 6}$.

\subsection{Equidistribution of root numbers}
Denote by $\epsilon_{r,s,t;\delta}$ the root number of the Jacobian variety $J_{r,s,t;\delta}$, and we know $\epsilon_{r,s,t;\delta}=\epsilon(\pi_{r,s,t;\delta})$. When $p=3$ ($r=s=t=1$), L. Mai \cite{Mai93} proved that the set $\{\delta \textrm{ cube-free}: \epsilon_{1,1,1;\delta}=+1 \}$ has density $1/2$ in the set of cube-free integers. For a general odd prime $p$, the curves $C_{1,1;p-2;\delta}$, $C_{1,(p-1)/2,(p-1)/2;\delta}$ and $C_{1,p-2,1;\delta}$ are hyperelliptic curves and the sign of the functional equation of L-functions of these curves can be obtained by \cite{Stoll93}. Following the strategy of L. Mai,  T. Jedrzejak \cite{Tomasz08} generalized the result of L. Mai for the families  $J_{1,1;p-2;\delta}$, $J_{1,(p-1)/2,(p-1)/2;\delta}$ and $J_{1,p-2,1;\delta}$ with $p\nmid \delta$. The constrain $p\nmid \delta$ comes from the calculations of root numbers in \cite{Stoll93}.

As an application of the general formulae for root numbers of $J_{r,s,t;\delta}$ given in Theorem \ref{root-number}, following the strategy and results in \cite{Mai93, Tomasz08}, we obtain the following density of root numbers of the families $J_{r,s,t;\delta}$ with $\delta$ varying in the set of  $p$-th power-free integers. 

Let $\BN^{(p)}$ be the set of $p$-th power-free integer which has positive density in the set of natural numbers (see for example \cite[Lemma 4]{Tomasz08}). Denote by $\BN^{(p)}_0\subset \BN^{(p)}$ the subset consisting of integers with $p$-adic valuation $0$. Then
\[\BN^{(p)}=\bigcup_{i=0}^{p-1}p^i\BN^{(p)}_0.\]
\begin{thm}\label{density}
\[\lim \limits_{X\ra +\infty} \frac{\#\{\delta\in \BN^{(p)}:\,\delta\leq X\textrm{ and } \epsilon_{r,s,t;\delta}=+1\}}{\#\{\delta\in \BN^{(p)}: \,\delta\leq X \}}=\frac{1}{2}.\]
\end{thm}

For $\delta\in \BN^{(p)}$, set 
\[\tau_p(\delta)=\#\left\{\ell\mid \delta:\,\ell\neq p \textrm{ and }\left(\frac{\ell}{p}\right)=-1\right\}.\]
If $\delta=p^i\delta_0$ with $(p,\delta_0)=1$, then it is clear that $\tau_p(\delta)$ is determined by $\delta_0$. Define $\alpha_p(\delta)=0,1$ by the formula
\[\epsilon_{r,s,t;\delta}=(-1)^{\alpha_p(\delta)+\tau_p(\delta)}.\]
\begin{lem}\label{dependence}
Suppose $\delta=p^i \delta_0\in \BN^{(p)}$ is a $p$-th power-free integer with $(p,\delta_0)=1$. Then $\alpha_p(\delta)$ depends on $(\ord_p(\delta), \delta_0\mod p^2)$.
\end{lem}
\begin{proof}
By the formulae for the root number $\epsilon_{r,s,t;\delta}$ in Theorem \ref{root-number}, if we write
\[\epsilon_{r,s,t;\delta}=(-1)^{\alpha_p(\delta)+\tau_p(\delta)},\]
then
\[(-1)^{\alpha_p(\delta)}=\left\{\begin{aligned}-\left(\frac{-rs\cdot \ord_p(\delta)}{p}\right),& \quad \textrm{ if } \ord_p(\delta)>0;\\
-\left(\frac{-rst\left(\frac{(r^rs^s(r+s)^t\delta^{r+s})^{p-1}-1}{p}\right)}{p}\right),&\quad \textrm{ if } u(r^rs^s(t-p)^t\delta^{r+s})=1;\\
\left(\frac{2}{p}\right),&\quad \textrm{ if }u(r^rs^s(t-p)^t\delta^{r+s})\geq 2.\end{aligned}\right.\]
It follows from these formulae that $\alpha_p(\delta)$ depends on $(\ord_p(\delta),\delta_0\mod p^2)$.
\end{proof}

\begin{proof}[Proof of Theorem \ref{density}]
In fact, we will prove for each $i=0,1,\cdots,p-1$ that 
\begin{eqnarray*}
&&\lim \limits_{X\ra +\infty} \frac{\#\{\delta\in p^i\BN^{(p)}_0:\,\delta\leq X \textrm{ and } \epsilon_{r,s,t;\delta}=+1\}}{\#\{\delta\in p^i\BN^{(p)}_0:\delta \leq X\}}\\
&=&\lim \limits_{X\ra +\infty} \frac{\#\{\delta\in \BN^{(i)}_0:\,\delta\leq X\textrm{ and } \epsilon_{r,s,t;p^i\delta}=+1\}}{\#\{\delta\in \BN^{(p)}_0:\delta \leq X\}}
=\frac{1}{2}.
\end{eqnarray*}
For fixed $i$, by Lemma \ref{dependence}, we can define $S_i^+$ resp. $S_i^-$ to be the subset of $(\BZ/p^2\BZ)^\times$ such that $\alpha_p(p^i\delta_0)=0,\textrm{ resp. }1$ for $\delta_0\in S_i^+$ resp. $S_i^-$. Then $S_i^+\cap S_i^-=\emptyset$ and $S_i^+\cup S_i^-=(\BZ/p^2\BZ)^\times$. For a sufficiently large integer $X$, denote
\[\BN^{(p)}_{0,X}=\{\delta\in \BN^{(p)}_{0}:\,\delta\leq X\}.\]
Then we have
\[\sum_{\begin{subarray}{c}\delta\in \BN^{(p)}_{0,X}\\\epsilon_{r,s,t;p^i\delta}=1\end{subarray}}1=\sum_{\begin{subarray}{c}\delta\in \BN^{(p)}_{0,X}\\\tau_p(\delta)\textrm{ even}\\ \delta\mod p^2\in S_i^+ \end{subarray}}1+\sum_{\begin{subarray}{c}\delta\in \BN^{(p)}_{0,X}\\\tau_p(\delta)\textrm{ odd}\\ \delta\mod p^2\in S_i^- \end{subarray}}1.\]
Let $\wh{(\BZ/p^2\BZ)^\times}$ be the group of characters of  the group $(\BZ/p^2\BZ)^\times$. Fix $k\in (\BZ/p^2\BZ)^\times$ and we have
\begin{eqnarray*}
\sum_{\begin{subarray}{c}\delta\in \BN^{(p)}_{0,X}\\(-1)^{\tau_p(\delta)}=\pm 1\\\delta\equiv k\mod p^2\end{subarray}}1&=&\sum_{\begin{subarray}{c}\delta\in \BN^{(p)}_{0,X}\\(-1)^{\tau_p(\delta)}=\pm 1\end{subarray}}\frac{1}{|\wh{(\BZ/p^2\BZ)^\times}|}\sum_{\chi\in \wh{(\BZ/p^2\BZ)^\times}}\chi(\delta)\chi^{-1}(k)\\
&=&\frac{1}{p^2-p}\sum_{\chi\in \wh{(\BZ/p^2\BZ)^\times}} \chi^{-1}(k)\sum_{\begin{subarray}{c}\delta\in \BN^{(p)}_{0,X}\\ (-1)^{\tau_p(\delta)}=\pm 1\end{subarray}} \chi(\delta)\\
&=&\frac{1}{p^2-p}\sum_{\chi\in \wh{(\BZ/p^2\BZ)^\times}} \chi^{-1}(k)\sum_{\begin{subarray}{c}\delta\in \BN^{(p)}_{0,X}\end{subarray}} \chi(\delta)\left(\frac{1\pm (-1)^{\tau_p(\delta)}}{2}\right)\\
&=&\frac{1}{2(p^2-p)}\sum_{\chi\in \wh{(\BZ/p^2\BZ)^\times}} \chi^{-1}(k)\sum_{\begin{subarray}{c}\delta\in \BN^{(p)}_{0,X}\end{subarray}} \chi(\delta)+O(\sqrt{X}\log^{p-2} X)\\
&=&\frac{1}{2(p^2-p)}\sum_{\begin{subarray}{c}\delta\in \BN^{(p)}_{0,X}\end{subarray}} \sum_{\chi\in \wh{(\BZ/p^2\BZ)^\times}} \chi^{-1}(k)\chi(\delta)+O(\sqrt{X}\log^{p-2} X)\\
&=&\frac{1}{2}\sum_{\begin{subarray}{c}\delta\in \BN^{(p)}_{0,X}\\ \delta\equiv k\mod p^2\end{subarray}}1+O(\sqrt{X}\log^{p-2} X).
\end{eqnarray*}
Here the $4$-th identify follows from \cite[Lemma 3]{Tomasz08}. Then by \cite[Lemma 5]{Tomasz08}, we obtain
\begin{eqnarray*}
\sum_{\begin{subarray}{c}\delta\in \BN^{(p)}_{0,X}\\\epsilon_{r,s,t;p^i\delta}=1\end{subarray}}1&=&\sum_{\begin{subarray}{c}\delta\in \BN^{(p)}_{0,X}\\\tau_p(\delta)\textrm{ even}\\ \delta\mod p^2\in S_i^+ \end{subarray}}1+\sum_{\begin{subarray}{c}\delta\in \BN^{(p)}_{0,X}\\\tau_p(\delta)\textrm{ odd}\\ \delta\mod p^2\in S_i^- \end{subarray}}1\\
&=&\frac{1}{2}\sum_{\begin{subarray}{c}\delta\in \BN^{(p)}_{0,X}\\ \delta\mod p^2\in S_i^+\end{subarray}}1+\frac{1}{2}\sum_{\begin{subarray}{c}\delta\in \BN^{(p)}_{0,X}\\ \delta\mod p^2\in S_i^-\end{subarray}}1+O(\sqrt{X}\log^{p-2} X)\\
&=&\frac{1}{2}\sum_{\begin{subarray}{c}\delta\in \BN^{(p)}_{0,X}\end{subarray}}1+O(\sqrt{X}\log^{p-2}X ).
\end{eqnarray*}

\end{proof}

\section{Selmer Groups}
Suppose $\delta \in \BQ^\times/\BQ^{\times p}$ is an integer. Recall for fixed $r,s,t$, we simply write $C_\delta$ for the curve
\[C_{r,s,t;\delta}:y^p=x^r(\delta-x)^s,\]
and write $J_\delta$ for the Jacobian variety of $C_\delta$. 

From the short exact sequence of $G_K$-modules
\[\xymatrix@C=11pt{0\ar[r]&J_{\delta}[\Pi]\ar[rr]&&J_{\delta}(\ov{\BQ})\ar[rr]^{\Pi}&&J_{\delta}(\ov{\BQ})\ar[r]&0},\]
taking Galois cohomologies, we have the Kummer map
\[\xymatrix{J_\delta(K)/\Pi J_\delta(K)\ar@{^(->}[r]^{\ \ \kappa}&\RH^1(K,J_\delta[\Pi]).}\]
For each place $V$ of $K$, similarly we have the local Kummer map
\[\xymatrix{J_\delta(K_V)/\Pi J_\delta(K_V)\ar@{^(->}[r]^{\ \ \kappa_V}&\RH^1(K_V,J_\delta[\Pi]).}\]
The Selmer group $\Sel({J_{\delta}}_{/K},\Pi)$ and the Shafarevich-Tate group $\Sha\left({J_\delta}_{/K}\right)$ are defined as
\[\Sel({J_{\delta}}_{/K},\Pi)=\{\xi\in \RH^1(K,J_\delta[\Pi]):\ \xi_V\in \Im(\kappa_V) \textrm{ for each } V.\}\]
\[\Sha\left({J_\delta}_{/K}\right)=\ker(\RH^1(K,J_\delta)\ra\prod_V\RH^1(K_V,J_\delta)).\]
It is well known that there is an exact sequence
\[0\ra J_\delta(K)/\Pi J_\delta(K)\xrightarrow{\kappa} \Sel({J_{\delta}}_{/K},\Pi)\ra \Sha\left({J_\delta}_{/K}\right)[\Pi]\ra 0.\]
We shall bound the Selmer group $\Sel({J_{\delta}}_{/K},\Pi)$, and necessarily we shall compute the local Kummer images $\Im(\kappa_V)$ for each place $V$ of $K$.

\subsection{Descent maps via Weil's idea}
Let $L\supset \mu_p$ be a field and $G_L$ its absolute Galois group. Let $\Div^0(C_\delta)(L)$ be the group of $L$-rational divisor of degree zero on the curve $C_\delta$. The group $J_\delta(L)$ is identified with the class group of $\Div^0(C_\delta)(L)$ modulo linear equivalence. The rational function $x$ on the curve  $C_\delta$ has divisor
\[\div(x)=p(0,0)-p(\infty)\]
where $(\infty)$ denotes the point on $C_\delta$ at infinity. To each divisor in $\Div^0(C_\delta)(L)$ away from the $ (0,0)$ and the infinity point $(\infty)$, we can associate an invariant as follows
\[\delta_L\left(\sum a_i P_i\right)=x\left(\sum a_i P_i\right)=\prod x(P_i)^{a_i}.\]
By Weil's reciprocity (see for example \cite{Coleman89}), on the principal divisors we find 
\[\delta_L(\div(f))=\left(\frac{f(0,0)}{f(\infty)}\right)^p\in L^{\times p}.\]
Hence $\delta_L$ induces a homomorphism, the invariant map,
\[\delta_L:J_\delta(L)\lra L^\times/L^{\times p}.\]
The kernel of this homomorphism turns out to be the subgroup generated by the classes $[1-\omega] D$ where $D\in \Div^0(C_\delta)(L)$, and hence we have an injective map
\[\delta_L:J_\delta(L)/\Pi J_\delta(L)\hookrightarrow L^\times/L^{\times p}.\]

From the short exact sequence of $G_L$-modules
\[\xymatrix@C=11pt{0\ar[r]&J_{\delta}[p]\ar[rr]&&J_{\delta}(\ov{L})\ar[rr]^{p}&&J_{\delta}(\ov{L})\ar[r]&0},\]
taking Galois cohomologies, we have the Kummer map
\[\xymatrix{J_\delta(L)/p J_\delta(L)\ar@{^(->}[r]^{\kappa_{L,p}}&\RH^1(L,J_\delta[p]).}\]
Denote
\[D_0=[(0,0)-(\infty)]\in J_\delta[\Pi].\]
Composed with the Weil pairing $e_p(\cdot,D_0)$, we obtain a map
\[e:\RH^1(L,J_\delta[p])\xrightarrow{e(\cdot,D_0)}\RH^1(L,\mu_p),\quad \xi\mapsto(\sigma\mapsto e_p(\xi(\sigma),D_0)).\]
By Hilbert Satz 90, we have an isomorphism
\[i:\RH^1(L,\mu_p)\simeq L^\times/L^{\times p}.\]
Finally, we get the descent map via the $p$-torsion point $D_0$:
\[i\circ e\circ \kappa_{L,p}:J_\delta(L)/p J_\delta(L)\ra L^\times/L^{\times p}\]
\begin{prop}\label{Weil}
The composition $i\circ e\circ \kappa_{L,p}:J_\delta(L)/p J_\delta(L)\ra L^\times/L^{\times p}$ factors through the invariant map $\delta_L:J_\delta(L)/\Pi J_\delta(L)\hookrightarrow L^\times/L^{\times p}$.
\end{prop}
\begin{proof}
Note any point of $J_\delta$ can be represented by a degree-$0$ divisor of $C_\delta$ which is not supported on $D_0$. Let $P$ and $Q$ be  divisors of degree $0$ on $C_\delta$ which are not supported on $D_0$. Suppose $P$ is $L$-rational and $pQ$ is linearly equivalent to $P$. Then the class $[P]$ represents a point in $J_\delta(L)$ and $[P]=p[Q]$.  If $[D_1]$ and $[D_2]$ are two torsion points of $J_\delta$ of order $p$, and suppose 
\[\div(f_1)=pD_1\textrm{ and }\div(f_2)=pD_2,\]
then the Weil pairing (see \cite[p. 284]{ACGH85})
\[e_p([D_1],[D_2])=\frac{f_2(D_1)}{f_1(D_2)}.\]
Let $\div(g)=pQ-P$. For any $\sigma\in G_L$,  $\div(g^\sigma )=pQ^\sigma-P$ and so $\div(g^\sigma/g)=pQ^\sigma-pQ.$
Therefore
\begin{eqnarray*}
e\circ \kappa_{L,p}([P])(\sigma)&=&e_p([Q^\sigma-Q],D_0)\\
&=&\frac{x(Q^\sigma-Q)}{g^\sigma/g(D_0)}\\
&=&\frac{\alpha^\sigma}{\alpha},
\end{eqnarray*}
where $\alpha=x(Q)/g(D_0)$. So we have
\[i\circ e\circ\kappa_{L,p}([P])\equiv \alpha^p\equiv \frac{x(pQ)}{g(pD_0)} \equiv x(P)\mod L^{\times p}.\]
\end{proof}

As the Jacobian variety of a curve, $J_\delta$ is principally polarized. We fix such a principal polarization to identify $J_\delta$ with its dual $J_\delta^\vee$.  Since the Rosati involution (w.r.t. the fixed polarization) is positive-definite, it coincides with the complex conjugation on the CM field $K$. Then the complex multiplication $[\ov{\Pi}]:J_\delta\ra J_\delta$ is the dual isogeny of $[\Pi]$, and $J_\delta[\Pi]=J_\delta[\ov{\Pi}]$.

Let $q=\left[{p}/{\ov{\Pi}}\right]: J_\delta[p]\ra J_\delta[\ov{\Pi}]$ be the natural quotient map. Then for any $x\in J_\delta[\Pi]$ and $y\in J_\delta[p]$, we have the following identity of Weil pairings:
\begin{equation}\label{WP}
e_\Pi(x,q(y))=e_p(x, y).
\end{equation}
By Proposition \ref{Weil}  we have the following commutative diagram:
\[\xymatrix{J_\delta(L)/\Pi J_\delta(L)\ar[rrrrd]^{\delta_L}&&&&\\
J_\delta(L)/pJ_\delta(L)\ar[r]^{\kappa_{L,p}}\ar[d]\ar[u]&\RH^1(L,J_\delta[p])\ar[r]\ar[d]^{q}\ar[r]^q&\RH^1(L,J_\delta[\ov{\Pi}])\ar[r]^{e'}&\RH^1(L,\mu_p)\ar[r]_i&L^\times/L^{\times p}\\
J_\delta(L)/\Pi J_\delta(L)\ar[r]^{\kappa_{L,\Pi}}&\RH^1(L,J_\delta(\Pi))\ar[ru]_{\mathrm{id}}&&&}\]
where $\kappa_{L,p}, \kappa_{L,\Pi}$ are Kummer maps associated to the isogenies $[p], [\Pi]$ respectively, $q$ is induced by the natural quotient map, $\id$ is the identity map,  and $e'$ is induced by composing the Weil pairing $e_\Pi(D_0,\cdot)$. Here we note $e=e'\circ q$ by (\ref{WP}).

Now we make the following identification
\begin{equation}\label{id}
i\circ e': \RH^1(L,J_\delta[\Pi])=L^\times/L^{\times p}.
\end{equation}
\begin{coro}
Under the identification (\ref{id}) the Kummer map $\kappa_{L,\Pi}$ coincides with the invariant map $\delta_L$.
\end{coro}

Hence in order to compute the images of the Kummer maps $\kappa_{L,\Pi}$, it suffices to compute the ranges of values of $\delta_L$, which is the group generated by norms of abscissas of points on the curve $C_\delta$ with coordinates in all possible algebraic extensions of $L$. However, by the following lemma,  it suffices to restrict oneself to extensions whose degree does not exceed the genus $g=(p-1)/2$ of $C_\delta$. 

\begin{lem}\label{bounding-field}
Let $C_{/L}$ be a smooth algebraic curve of genus $g$ with an $L$-rational point $P_0$. In each class of $L$-rational divisors of degree zero, there is a divisor of the form $P_1+P_2+\cdots+P_g-gP_0$, and consequently, the degree of the field, generated by the coordinates of each $P_i$, does not exceed $g$.
\end{lem}
\begin{proof}
Since $P_1+P_2+\cdots+P_g-gP_0$ is $L$-rational, for each $1\leq i\leq g$, all the Galois conjugates of $P_i$ occur in  $P_1,P_2,\cdots,P_g$, and hence the field generated by the coordinates of each $P_i$ has degree $\leq g$. Now it suffices to prove there exists a divisor of such form. Let $D$ be an arbitrary $L$-rational point of degree zero, and let $P_0$ be an arbitrary point in $C(L)\neq \emptyset$. By Riemmann-Roch thoerem,
\[l(D+g P_0)\geq \deg(D+g P_0)-g+1= 1.\]
Therefore $D+gP_0$ is linearly equivalent to an effective divisor of degree $g$, and the lemma follows.
\end{proof}


Before we go further, we establish a dimension formula for the local Kummer images, which will save us much energy when compared with explicit calculations in \cite{Fad61}. Let $A_{/L}$ be an abelian variety over a nonarchimedean local field $L$ of characteristic $0$. Let $\alpha:A\ra A$ be an isogeny over $L$ of degree $n$, and denote by  $\alpha^\vee:A^\vee\ra A^\vee$ the dual isogeny. We know the Cartier duality  $A[\alpha]^\vee=A^\vee[\alpha^\vee]$. By the local Tate duality for finite groups \cite[Theorem 8.1]{Tate-Galois} \cite[Corollary I-3.2,  III-6.10]{Milne-ADT06}, we have a perfect pairing
\[\RH^1(L,A[\alpha])\times \RH^1(L,A^\vee[\alpha^\vee])\lra \mu_n.\]
Under the Kummer maps, we view $A(L)/\alpha A(L)$ and $A^\vee(L)/\alpha^\vee A^\vee(L)$ as subgroups of $\RH^1(L,A[\alpha])$ and $\RH^1(L,A^\vee[\alpha^\vee])$ respectively.

\begin{lem}\label{complements}
The subgroups  $A(L)/\alpha A(L)$ and $A^\vee(L)/\alpha^\vee A^\vee(L)$ are exact annihilators under the local Tate pairing.
\end{lem}
\begin{proof}
Using Weil pairing, we have a diagram
\[
\xymatrix{0\ar[r]&A(L)/\alpha A(L)\ar[d]\ar[r]&\RH^1(L,A[\alpha])\ar[r]\ar[d]&\RH^1(L,A)[\alpha]\ar[r]\ar[d]&0\\
0\ar[r]&\RH^1(L,A^\vee)[\alpha^\vee]^*\ar[r]&\RH^1(L,A^\vee[\alpha^\vee])^*\ar[r]&A^\vee(L)/\alpha^\vee A^\vee(L)^*\ar[r]&0}
\]
with exact rows, where $^*$ denotes the Pontrjagin dual. By the local Tate dualities for finite groups and abelian varieties \cite[Theorem I-2.1, 3.2,  III-6.10, Corollay I-2.3, 3.4]{Milne-ADT06},  all the vertical maps are isomorphisms. Hence,  $A(L)/\alpha A(L)$ and $A^\vee(L)/\alpha^\vee A^\vee(L)$ are exact annihilators.
\end{proof}

\begin{prop}\label{dim}
We have
\[\dim_{\BF_p} \Im(\kappa_{L,\Pi})=\frac{\dim_{\BF_p}\RH^1(L,J_\delta[\Pi])}{2}.\]
\end{prop}
\begin{proof}
Take $A=A^\vee=J_\delta$ and $\alpha=\Pi, \alpha^\vee=\ov{\Pi}$ in Lemma \ref{complements}. We have the local Tate duality
\begin{equation}
\label{pairing2}\RH^1(L,J_\delta[\Pi])\times \RH^1(L,J_\delta[\ov{\Pi}])\lra \BQ/\BZ.
\end{equation}

By Lemma \ref{complements}, the Kummer images of $J_\delta(L)/\Pi J_\delta(L)$  and $J_\delta(L)/\ov{\Pi} J_\delta(L)$ are orthogonal complements of each other.  Since the unit $\ov{\Pi}/\Pi\in \CO_K^\times$ induces an isomorphism $[\ov{\Pi}/\Pi]:J_\delta(L)/\Pi J_\delta(L)\ra J_\delta(L)/\ov{\Pi} J_\delta(L)$ by complex multiplication, by the following commutative diagram
\[\xymatrix{
J_\delta(L)/\Pi J_\delta(L)\ar[rd]^{\kappa_{L,\Pi}}\ar[dd]^{[\ov{\Pi}/\Pi]}&\\
&\RH^1(L,J[\Pi])\\
J_\delta(L)/\ov{\Pi}J_\delta(L)\ar[ru]^{\kappa_{L,\ov{\Pi}}}&\\}\]
we see $\Im(\kappa_{L,\Pi})$ and $\Im(\kappa_{L,\ov{\Pi}})$ coincide. Hence under the nondegenerate pairing (\ref{pairing2}), $\Im(\kappa_{L,\Pi})$ is the exact annihilator of itself, and the proposition follows.
\end{proof}

For each  place $V$ of $K$, denote $\delta_{K_V}$ simply by $\delta_V$. The Kummer images at the infinite places are trivial and we shall compute the images $\Im(\delta_V)$ for all finite places $V$,  which is the group generated by norms of abscissas of points on the curve $C_\delta$ with coordinates in all possible algebraic extensions of $K_V$. However, by Riemann-Roch,  it suffices to restrict oneself to extensions whose degree does not exceed the genus $g=(p-1)/2$ of $C_\delta$.


If 
\[(r',s',t')=h(r,s,t)+p(i,j,k),\]
with $r',s',t'>0,$ and $r'+s'+t'=p$. Then we have a commutative diagram
\[\xymatrix{&F_\delta\ar[ld]_{\varphi_{r,s,t}}\ar[rd]^{\varphi_{r',s',t'}}&\\C_{r,s,t;\delta}\ar[rr]&&C_{r',s',t';\delta}\ ,}\]
where the horizontal arrow is
\[(x,y)\mapsto(x,y^hx^i(\delta-x)^j).\]
Since $h$ is prime to $p$, this map is a birational equivalence. 
Without loss of generality, we may fix $s=1$ and take $1\leq r\leq p-2$.  We are reduced to investigate the norms from $L$ to $K_V$ of the abscissas of $L$-points on the curves
\begin{equation}\label{eqn-1}
y^p=x^r(\delta-x), \quad 1\leq r\leq p-2,
\end{equation}
where $[L:K_V]\leq (p-1)/2$.

For our convenience, we denote $K_V$ by $k$, $L$--an algebraic extension of the field $k$ of degree $[L:k]\leq (p-1)/2$. The uniformizers of $k$ and $L$ are denoted by $\lambda$ and $\Lambda$. The valuations of elements $x$ in $L$ will be denoted by $\ord(x)$, such that $\ord(\lambda)=1$. Thus, valuations can be fractional numbers, with denominators not exceeding $(p-1)/2$.

\subsection{Local descents outside $p$}
Let $\ell\neq p$ be a prime number, and let $V\mid \ell$ be a place of $K$. Since the extension $K/\BQ$ is unramified at $\ell$, we may take the uniformizer $\varPi_V=\ell$  for $K_V$.
\[K_V^\times/K_V^{\times p}=\langle \ell\rangle^{\BZ/p\BZ}\times \BF_{V}^\times/\BF_V^{\times p}\]
has dimension $2$ over $\BF_p$. Then by Proposition \ref{dim},
\begin{equation}\label{selmer-dimension}
\dim_{\BF_p}\Im(\delta_V)=1.
\end{equation}

\begin{thm}\label{desc-off-p}
Suppose $(r,s,t)=(r,1,p-r-1)$ and $\delta\in \BQ^\times/\BQ^{\times p}$ is an integer. Let $\ell$ be a prime number not equal to $p$, and let $V\mid \ell$ be a place of $K$. If $\ell\nmid \delta$, then 
\[\Im(\delta_V)=\BF_V^\times/\BF_V^{\times p}.\]
If $\ell\mid \delta$, then 
\[\Im(\delta_V)=\langle \delta \rangle^{\BZ/p\BZ}.\]
\end{thm}
\begin{proof}
If $\ell\nmid \delta$, then $J_\delta$ has good reduction at $V$. By \cite[Lemma 6]{GP12}, we have
\[\Im(\kappa_V)=\RH^1_\ur(K_V,J_\delta[\Pi]),\]
and hence 
\[\Im(\delta_V)=\BF_V^\times/\BF_V^{\times p}.\]

Next we assume $\ell\mid \delta$. we shall investigate the norms of abscissas of points on the curve $C_\delta$, and we present two cases:

\noindent{\bf Case I.} $ \ord(x)<\ord(\delta).$ Substituting by
\[x_1=\frac{\delta}{x}，\quad y_1=-\frac{y}{x}\]
in Equation (\ref{eqn-1}), we get
\[x_1^{r_1}(1-x_1)=y_1^p\delta^{r_1},\quad r_1=p-r-1.\]
Let $1\leq r_2\leq p-2$ be such that $r_1r_2=1+pk$ for some $k$.  Putting
\[x_2=x_1,\quad y_2=y_1^{r_2}\delta^kx_1^{-k},\]
we have
\[x_2(1-x_2)^{r_2}=y_2^p\delta.\]
Then
\[x_2-r_2x_2^2+\frac{r_2(r_2-1)}{2}x_2^3+\cdots=y_2^p\delta.\]
Since $\ord(x_2)>0$ and, consequently, $\ord(y_2^p\delta)>0$, by \cite[Lemma 1]{Fad61},
\[x_2=y_2^p\delta(1+r_2y_2^p\delta+\cdots).\]
Since $\ord(y_2^p\delta)>0$,
\[\N_{L/K_V}(1+r_2y_2^p\delta+\cdots)\in 1+\ell\CO_{K,V}\subset K_V^{\times p}.\]
Then 
\[\N_{L/K_V}(x)=\delta^{[L:K_V]}\N_{L/K_V}(x_2)^{-1}\equiv 1\mod K_V^{\times p}.\]

\noindent{\bf Case II.} $ \ord(x)\geq\ord(\delta).$ Put 
\[x_1=\frac{x}{\delta},\quad y_1=y.\]
Then Equation (\ref{eqn-1}) transforms into
\[x_1^r(1-x_1)=y_1^p\delta^{-(r+1)}.\]
First we claim that either $\ord(x_1)>0$ or $\ord(1-x_1)>0$.  On the contrary, if $\ord(x_1)=\ord(1-x_1)=0$, then
\[p\cdot \ord(y_1)=(r+1)\ord(\delta).\]
By our choice of $\delta$, we have $1\leq \ord(\delta)\leq p-1$. Since the denominator of $\ord(y_1)$ doesn't exceed $(p-1)/2$, we conclude that $\ord(y_1)$ is an integer and $p\mid \ord(\delta)$, which contradicts with our choice of $\delta$.

If $\ord(1-x_1)>0$, then 
\[\N_{L/K_V}(x_1)=\N_{L/K_V}(1+(x_1-1))\in 1+\ell \CO_{K,V}\subset K_V^{\times p}.\]
In this case,
\[\N_{L/K_V}(x)=\delta^{[L:K_V]}\N(L/k)(x_1)\equiv \delta^{[L:K_V]} \mod K_V^{\times p}.\]

Next suppose $\ord(x_1)>0$. Let $1\leq r_1\leq p-2$ be such that $rr_1=1+pk$ for some $k$ and put
\[x_2=x_1,\quad y_2=\frac{y_1^{r_1}}{(\delta x_1)^{k}}.\] 
Then we have
\[x_2(1-x_2)^{r_1}=y_2^{p}\delta^{-(r_1+1)},\]
and hence
\[x_2-r_1x_2^2+\cdots=y_2^p\delta^{-(r_1+1)}.\]
Since $\ord(x_2)=\ord(x_1)>0$, and consequently $\ord(y_2^p\delta^{-(r_1+1)})>0$, by \cite[Lemma 1]{Fad61},
\[x_2=y_2^p\delta^{-(r_1+1)}(1+r_1y_2^p\delta^{-(r_1+1)}+\cdots)\]
In this case,
\[\N_{L/K_V}(x)=\delta^{[L:K_V]}\N_{L/K_V}(x_2)\equiv\delta^{-r_1[L:K_V]}\mod K_V^{\times p}.\]
Therefore, we see that , in any case, the norm $\N_{L/K_V}(x)$ falls into the subgroup $\langle \delta\rangle \subset K_V^\times/K_V^{\times p}$. From (\ref{selmer-dimension}), $\Im(\delta_V)$ has dimension $1$ over $\BF_p$, we conclude that $\Im(\delta_V)=\langle \delta \rangle\subset K_V^\times/K_V^{\times p}.$ 

In the following, we explicitly construct a divisor whose image under $\delta_V$ is a  non-trivial power  of $\delta$ in $k^\times/k^{\times p}$. First we assume $\ell\nmid (r+1)$ and construct points with $\ord(x)=0$ which gives a point in Case I.
Now take $x_0\in \BF_\ell^{\times p}$ and put $x=x_0+z$. Then the equation determining $z$ is 
\[x_0^{r-1}(r\delta-(r+1) x_0)z+\cdots =y^p-x_0^{r}(\delta-x_0).\]
Note  $x_0^{r-1}(r\delta-(r+1) x_0)$ is an $\ell$-adic unit in $\BQ_\ell^\times$.
By \cite[Lemma 1]{Fad61}, the solvability of $z$ with $\ord(z)>0$ is equivalent to
\[y^p\equiv x_0^r(\delta-x_0)\equiv -x_0^{r+1}\mod \Lambda.\]
By our choice of $x_0\in \BF_\ell^{\times p}$, we can take $y_1\in \BF_\ell^\times$ such that
\[y_1^p=-x_0^{r+1}.\]
Then the equation 
\[x_0^{r-1}(r\delta-(r+1) x_0)z+\cdots =y_1^p-x_0^{r}(\delta-x_0)\]
has a unique solution
\[z_1=\frac{y_1^p-x_0^r(\delta-x_0)}{x_0^{r-1}(r\delta-(r+1)x_0)}+b_2\left(\frac{y_1^p-x_0^r(\delta-x_0)}{x_0^{r-1}(r\delta-(r+1)x_0)}\right)^2+\cdots \in \BQ_\ell,\]
with $\ord(z_1)>0$. If we put $x_1=x_0+z_1$, then $(x_1,y_1)$ is an $\BQ_\ell^\times$-point on the curve 
\[x^r(\delta-x)=y^p.\]
Now
\[\N_{L/k}(x_1)=x_1\equiv x_0\equiv 1\mod k^{\times p}.\]
On the other hand, $(\delta,0)$ also gives a $\BQ_\ell$-point which gives a point in Case II, and 
\[\delta_V([(\delta,0)-(x_1,y_1)])=\delta \mod k^{\times p}.\]

For the primes $\ell\mid (r+1)$, we use the following trick. Let $r_1>0$ be such that 
\[rr_1=pk+1, \textrm{ and }\ell\nmid (r_1+1).\]  
The substitution $x_1=x$, $y_1=y^{r_1}x^{-k}$ realizes a birational map from the curve $x^r(\delta-x)=y^p$ to $x_1(\delta-x_1)^{r_1}=y_1^p$. 
Again, put $x_1=x_0+z$, and the equation determining $z$ is 
\[(\delta-x_0)^{r_1-1}(\delta-(r_1+1) x_0)z+\cdots =y_1^p-x_0(\delta-x_0)^{r_1}.\]
Taking $x_0\in \BF_\ell^{\times p}$, $(\delta-x_0)^{r_1-1}(\delta-(r_1+1) x_0)$ is an $\ell$-adic unit in $\BQ_\ell^\times$. By a same construction as above, we see there exists a rational divisor with image equal to $\delta \mod k^{\times p}$ under the Kummer map $\delta_V$.

\end{proof}

\subsection{Local descents at $p$}
In this subsection we suppose $p\geq 5$ and impose the following condition 
\[p\nmid \delta.\]
We  follow \cite{Fad61} with appropriate modifications.
It will be convenient to represent the field $K_\Pi$ as $\BQ_p(\lambda)$, where $\lambda^{p-1}+p=0$.
As generators of $K_\Pi^\times/K_\Pi^{\times p}$, we take 
\begin{eqnarray*}
u_0&=&\lambda,\\
u_1&=&\omega,\\
u_i&=&e^{\lambda^i}, \quad 2\leq i\leq p.
\end{eqnarray*}
Let $\mu:\Gal(K/\BQ)\ra \BZ_p^\times$ be the Teichm\"uller character:
\[\omega^{\sigma}=\omega^{\mu(\sigma)}.\]
Then $\mu$ has values in $\BF_p^\times$, and
\[\lambda^\sigma=\mu(\sigma)\lambda.\]
If $M$ is a $\BZ_p[\Gal(K/\BQ)]$-module, we denote by $M(i)$ the submodule of $M$ consisting of all $x\in M$ such that $x^\sigma=\mu^i(\sigma)x$ for all $\sigma\in \Gal(K/\BQ)$. 

Denote $U_i$ to be the one-dimensional subspace of $K_\Pi^\times/K_\Pi^{\times p}$ generated by $u_i$ for $0\leq i\leq p$. Then  $\Gal(K/\BQ)$ acts on $U_i$ as $\mu^i$, and 
\begin{eqnarray*}
K_\Pi^\times/K_\Pi^{\times p}(0)&=&U_0\times U_{p-1},\\
K_\Pi^\times/K_\Pi^{\times p}(1)&=&U_1\times U_p,\\
K_\Pi^\times/K_\Pi^{\times p}(i)&=&U_i,\quad 2\leq i\leq p-2.
\end{eqnarray*}

Recall we have denoted by $k=K_\Pi$, $L$--an algebraic extension of the field $k$ of degree $[L:k]\leq (p-1)/2$. Denote $k_1$--the maximal unramified subfield of $L$ over $\BQ_p$, $L_1=kk_1$--the maximal unramified subfield of $L$ over $k$. The uniformizers of $k$ and $L$ are denoted by $\lambda$ and $\Lambda$. The valuations of elements $x$ in $L$ will be denoted by $\ord(x)$, such that $\ord(\lambda)=1$. 

We shall investigate the norms from $L$ to $k$ of the abscissas of $L$-points on the curves
\[y^p=x^r(\delta-x), \quad 1\leq r\leq p-2.\]
Here exact five cases can be presented:
\begin{itemize}
\item[I.] $\ord(x)<0$;
\item[II.] $\ord(x)>0$;
\item[III.] $x\equiv \delta \mod \Lambda$;
\item[IV.] $\ord(x)=0$, and $x\not \equiv \delta, r\delta/(r+1)\mod \Lambda$;
\item[V.] $x\equiv r\delta/(r+1)\mod \Lambda$.
\end{itemize}

\noindent{\bf  Case I.}
The substitutions $x_1=x^{-1}$ and $y_1=-x^{-1}y$ realize a birational map from the curve $y^p=x^r(\delta-x)$ onto the curve 
$y_1^p=x_1^{r_1}(1-\delta x_1)$, where $r_1=p-r-1$. The latter equation can be transformed into $y_2^p=x_2(1-\delta x_2)^{r_2}$ by the substitutions $x_2=x_1$ and $y_2=y_1^{r_2}x_2^{-k}$, where $r_1r_2=1+kp$, $1\leq r_2\leq p-2$. It is clear from the substitutions that the norms of the abscissas of the corresponding $L$-points on $y^p=x^r(\delta-x)$ and $y^p_2=x_2(1-\delta x_2)^{r_2}$ fall into one same subgroup of $k^\times/k^{\times p}$. Therefore we only need to consider $L$-points on 
\[y^p=x(1-\delta x)^r, \quad 1\leq r\leq p-2,\]
with $\ord(x)>0$. 

Expanding in powers of $x$, we obtain
\[x-r\delta x^2+\frac{r(r-1)}{2}\delta^2x^3+\cdots =y^p.\]
Since $\ord(y)>0$, by \cite[Lemma 1]{Fad61}, 
\[x=y^p(1+b_2 y^{p}+b_3y^{3p}+\cdots).\]
with $p$-adic integral rational coefficients $b_n$ and $b_2=r\delta$.
Then 
\[\N_{L/k}(x)\equiv \N_{L/k}(1+b_2y^p+b_3y^{2p}\cdots)\mod k^{\times p}.\]
By \cite[Lemmas 3, 4]{Fad61}, 
\[\Tr_{L/k}(y^p)\equiv \Tr_{L/k}(y^{2p})\equiv \cdots \equiv 0\mod \lambda^p.\] 
By \cite[Lemma 5]{Fad61}, the norm $\N_{L/k}(x)$ falls into the subgroup $U_p$. For $K=k$, $y=\lambda $, we have
\[\N_{L/k}(1+b_2y^p+b_3y^{2p}\cdots)=1+b_2y^p+b_3y^{2p}\cdots\equiv 1+r\delta \lambda^{p}\equiv u_p^{r\delta}\mod  \lambda^{2p}. \]
It follows that there exists $\N_{L/k}(x)$ which is a non-trivial element in  $U_p$.

\noindent{\bf Cases II and III.} \label{case1-23}
It is not necessary to investigate the cases II and III separately. The substitutions $x_1=\delta-x$ and  $y_1=y^{r_1}x^{-k}$ realize a birational map from the curve $y^p=x^r(\delta-x)$ on to the curve $y_1^p=x_1^{r_1}(\delta-x_1)$ where $rr_1=1+k p$, $1\leq r_1\leq p-2$.  Now Case II and Case III change places, and the norms of the abscissas of the corresponding points again fall into one same subgroup of $k^\times/k^{\times p}$, because
\[\N_{L/k}(x_1)=\N_{L/k}(\delta-x)=\N_{L/k}(y)^p\N_{L/k}(x)^{-r}.\]

Now suppose $x\equiv \delta\mod \Lambda$. Set $x=\delta+z$, $\ord(z)>0$. After simple transformations, we have
\[z+r\delta^{-1}z^2+\frac{r(r-1)}{2}\delta^{-2}z^3+\cdots=-y^p\delta^{-r}.\]
By \cite[Lemma 1]{Fad61}, 
\[z=-y^p\delta^{-r}+b_2y^{2p}\delta^{-2r}+\cdots.\]
Since $\ord(y)>0$, we have, by \cite[Lemmas 3, 4]{Fad61}, 
\[\Tr_{L/k}(y^p)\equiv \Tr_{L/k}(y^{2p})\equiv \cdots \equiv 0\mod \lambda^p.\] 
By \cite[Lemma 5]{Fad61}, the norm $N_{L/k}(x)$ falls into the subgroup generated by $u_{p-1}$ and $u_p$. For $L=k$, $y=\lambda $, we have
\[\N_{L/k}(x)=x\equiv \delta(1-\lambda^{p}\delta^{-r-1})\equiv \delta u_p^{-\delta^{-r-1}} \mod \lambda^{2p}. \]
Since $\delta$ is Galois invariant, it belongs to $K_\Pi^\times/K_\Pi^{\times p}(0)=U_0\times U_{p-1}$. Since $\delta$ is a $p$-adic unit, it is a power of $u_{p-1}$ in $K_\Pi^\times/K_\Pi^{\times p}$. Then there exists some $\N_{L/k}(x)$ which has nontrivial exponents at $u_{p-1}$ and $u_p$.

\noindent{\bf Case IV.} \label{case1-4}
Suppose $x\not \equiv \delta, r\delta/(r+1)\mod \Lambda$. Denote by $x_0$ the residue of $x \mod \Lambda$,  and put $x=x_0+z$.  It is known that we can take $x_0\in k_1$. Expanding in powers of $z$, we get
\[x_0^{r-1}(r\delta-(r+1)x_0)z+\cdots=y^p-x_0^r(\delta-x_0).\]
The coefficients on the left are integers in $k_1$. The resulting equation will have a solution for $\ord(z)>0$ only if 
\[y^p\equiv x_0^r(\delta-x_0)\mod \Lambda,\]
and, equivalently,
\[y\equiv (x_0^r(\delta-x_0))^{p^{f-1}}\mod \Lambda,\]
where $f$ is the degree of $k_1$ over $\BQ_p$. 

Set $y_1= (x_0^r(\delta-x_0))^{p^{f-1}}\in k_1$. Then there exists a unique solution $z_1\equiv 0\mod \Lambda$ of the equation
\[x_0^{r-1}(r\delta-(r+1)x_0)z+\cdots=y_1^p-x_0^r(\delta-x_0),\]
and $z_1\in k_1$. Then $x_1=x_0+z_1$ and $y_1$ for a solution of the equation $x^r(\delta-x)=y^p$ in $k_1$.

Now set $x=x_1+z$ and perform similar transformations:
\[x_1^{r-1}(r\delta-(r+1)x_0)z+\cdots=y^p-y_1^p.\]
Put
\[y=y_1(1+v),\quad A=\frac{y_1^p}{x_1^{r-1}(r\delta-(r+1)x_1)}.\] 
The equation changes into the form
\[z+c_2 z^2+\cdots=A((1+v)^p-1)=Av^p+p w,\]
\[v,w\in L,\quad \ord(v),\ord(w)>0.\]
Consequently,
\[z=Av^p+A_1 v^{2p}+\cdots+p w',\]
\[A,A_1,\cdots \in k_1, \ord(A),\ord(A_1)\cdots\geq0, w'\in L,\ord(w')>0.\]
We calculate $\Tr_{L/L_1}(z)$. Since $L_1$ is unramified over $k$, and $\ord(v),\ord(w')>0$,
\[\Tr_{L/L_1}(v)\equiv \Tr_{L/L_1}(w')\equiv 0\mod \lambda.\] 
Applying \cite[Lemmas 3, 4]{Fad61}, we get
\[\Tr_{L/L_1}(z)\equiv \Tr_{L/L_1}(z^2)\equiv \cdots \equiv 0\mod \lambda^p.\]
Therefore
\[\N_{L/L_1}(x)=\N_{L/L_1}(x_1)\N_{L/L_1}(1+x_1^{-1}z)\equiv x_1^{[L:L_1]}\mod \lambda^p.\]
Consequently, 
\[\N_{L/k}(x)\equiv \N_{L_1/k}(x)^{[L:L_1]}\mod \lambda^p.\]
But $\N_{L_1/k}(x_1)=\N_{k_1/\BQ_p}(x_1)\in \BQ_p$, and hence $\N_{L/k}(x)$ is congruent to rational numbers modulo $\lambda^p$. Consequently, $\N_{L/k}(x)$ falls into the subgroup of $k^\times/k^{\times p}$ generated by $u_{p-1}$ and $u_p$.

Next we show that for $p\geq 5$, $x$ can be taken so that $\N_{L/k}(x)$ is a non-trivial element in $U_{p-1}$. Let
\[\varphi(x)=\frac{(1-x)^p-1+x^p}{p}=x\left(x^{p-2}-\cdots+\frac{p-1}{2}x-1\right),\]
and define a map:
\[\varphi:\BF_p\lra \BF_p,\quad x\mapsto \varphi(x).\]

\begin{lem}\label{choice}
Suppose $p\geq 7$. Then the cardinalities $|\varphi^{-1}(0)|\leq p-3$ and $|\varphi^{-1}(a)|\leq p-4$ for all $a\neq 0$.
\end{lem}
\begin{proof}
It is known that $\varphi(x)$ is divisible by $x^2-x+1$ if $p\equiv 5\mod 6$, and $\varphi(x)$ is divisible by $(x^2-x+1)^2$ if $p\equiv 1\mod 6$. Therefore the number of distinct solutions of $\varphi(x)=0$ is at most $p-3$. If $p\equiv 1\mod 6$, then $x^2-x+1$ has two distinct roots which are different from $0,1$. Hence $|\varphi^{-1}(0)|\geq 4$, and the lemma holds in this case. In the following we give a uniform proof for prime $p>5$. The lemma is valid for $p=7$, and we suppose $p>7$. 

It is easily seen that 
\[\varphi(x)=\varphi(1-x),\quad x\in \BF_p.\]
By direct calculation, we have
\[\varphi(0)=\varphi(1)=0,\varphi(2)=\frac{2^p-2}{p},\textrm{ and }\varphi\left(\frac{p+1}{2}\right)=\frac{2^{-(p-1)}-1}{p}.\] 
And $\varphi(2)=\varphi((p+1)/2)$ if and only if 
\[\frac{2^p-2}{p}=0.\]

If $\varphi(2)=\varphi((p+1)/2)$, then $0,1,2,p-1, (p+1)/2$ are distinct elements in the fibre $\varphi^{-1}(0)$, and hence the cardinality of other fibers are $\leq p-5$.

Next suppose $\varphi(2)\neq \varphi((p+1)/2)$, and hence
\[\frac{2^p-2}{p}\neq 0, \textrm{ i.e. } \varphi(2)\neq 0.\]
Since  the $\varphi(0)=\varphi(1)$, $\varphi(2)=\varphi(p-1)$ and $\varphi((p+1)/2)$ are all distinct,  we see
\[|\varphi^{-1}(0)|\leq p-3, \textrm{ and } |\varphi^{-1}(a)|\leq p-4 \textrm{ for any } a\neq 0\textrm{ and } \varphi(2).\]
The inequality
\[|\varphi^{-1}(\varphi(2))|\leq p-5\]
follows from the claim
\[\varphi(2)\neq \varphi(3), \textrm{ or }\varphi(2)\neq\varphi(4).\]
By the assumption $p>7$, we have $4<(p-1)/2$. 
Suppose the contrary
\[\varphi(2)= \varphi(3)=\varphi(4),\]
equivalently,
\[\frac{2^p-2}{p}=\frac{-2^p-1+3^p}{p}=\frac{-3^p-1+4^p}{p}.\]
Then
\[3^p\equiv 2^{p+1}-1\mod p^2,\quad 3^p\equiv 2^{p-1}(2^p+1)\mod p^2.\]
Suppose 
\[2^{p-1}=1+pt,\quad t\in \BZ.\]
Then 
\[2^{p-1}(2^p+1)\equiv 3+5pt\mod p^2,\quad 2^{p+1}-1\equiv 3+4pt\mod p^2.\]
It follows that $t\equiv 0\mod p$. Hence
\[2^{p-1}-1\equiv 0\mod p^2\textrm{ and }\varphi(2)=0,\]
which contradicts with the assumption $\varphi(2)\neq \varphi((p+1)/2)$.

\end{proof}
We  choose $x_0\in \BQ_p$ so that $x_0^p=x_0$, and therefore $y_1=x_0^r(\delta-x_0)$. The equation for $z_1$ is 
\begin{eqnarray*}
x_0^{r-1}(r\delta-(r+1)x_0)z+\cdots &=&y_1^p-y_1\\
&=&x_0^{r}((\delta-x_0)^p-(\delta-x_0))\\
&=&px_0^r \left(\frac{\delta^p-\delta}{p}+\delta^p \varphi\left(\frac{x_0}{\delta}\right)\right).\end{eqnarray*}
Note that $\varphi(0)=\varphi(1)=0$. We can choose $x_0$ with additional properties:
\[\frac{x_0}{\delta}\not \equiv 0, 1, \frac{r}{r+1}\mod p,\]
and 
\[\frac{\delta^p-\delta}{p}+\delta^p\varphi\left(\frac{x_0}{\delta}\right)\not \equiv 0\mod p.\]
When $p\geq 7$, this is achieved by  Lemma \ref{choice}. Such a choice of $x_0$ holds for the case $p=5$ just by noting
\[\varphi(0)=\varphi(1)=0,\quad \varphi(2)=\varphi(4)=1,\quad \varphi(3)=2,\]
and
\[\frac{\delta^5-\delta}{5}+\delta^5\not \equiv 0\mod 5.\]

Then by \cite[Lemma 1]{Fad61},
\[z_1\equiv \frac{px_0 \left(\frac{\delta^p-\delta}{p}+\delta^p \varphi\left(\frac{x_0}{\delta}\right)\right)}{r\delta-(r+1)x_0}\mod p^2,\]
and 
\begin{eqnarray*}
x_1=x_0+z_1&\equiv& x_0\left(1+\frac{p \left(\frac{\delta^p-\delta}{p}+\delta^p \varphi\left(\frac{x_0}{\delta}\right)\right)}{r\delta-(r+1)x_0}\right)\mod p^2.
\end{eqnarray*}
Further, $x_0^p=x_0$, so that 
\[x_1\equiv u_{p-1}^{c} \mod k^{\times p},\]
where
\[c\equiv -\frac{\frac{\delta^p-\delta}{p}+\delta^p \varphi\left(\frac{x_0}{\delta}\right)}{r\delta-(r+1)x_0}\not\equiv 0\mod p.\]

\noindent{\bf Case V.}\label{case1-5}
In this case,
\[x\equiv \frac{r\delta}{r+1}\mod \Lambda.\]
We set $x=r\delta/(r+1)+z$, and expanding the left-hand side of $x^r(\delta-x)=y^p$ by the Taylor expansion, we obtain
\begin{equation}\label{eqn1}
-\frac{(r\delta)^{r-1}}{2(r+1)^{r-2}}z^2+\cdots=y^p-\frac{r^r\delta^{r+1}}{(r+1)^{r+1}},
\end{equation}
where the coefficients with respect to powers of $z$ are $p$-adic integral rational numbers. Equation (\ref{eqn1}) can have a solution for $z\equiv 0\mod \Lambda$ only when
\[y\equiv \frac{r^r\delta^{r+1}}{(r+1)^{r+1}}\mod \Lambda.\]
We set 
\[y=\frac{r^r\delta^{r+1}}{(r+1)^{r+1}}(1+v).\]
Then Equation (\ref{eqn1}) transforms into
\begin{equation}\label{eqn2}
z^2(1+a_3z+\cdots)=A\left((1+v)^p-1\right)-pB,
\end{equation}
where
\[A=-\left(\frac{r^r\delta^{r+1}}{(r+1)^{r+1}}\right)^p\cdot \frac{2(r+1)^{r-2}}{(r\delta)^{r-1}},\]
\[B=p^{-1}\left(\left(\frac{r^r\delta^{r+1}}{(r+1)^{r+1}}\right)^{p-1}-1\right)\cdot \frac{2r\delta^2}{(r+1)^3}.\]
The coefficients $A,B,a_3,\cdots$ are all $p$-adic integral rational numbers. Further investigations should be conducted by distinguishing several subcases.

\noindent{\bf V-A. $\ord(v)<1$.} Then 
\[\ord(v)\leq 1-\frac{2}{p-1},\textrm{ and }\ord(v^p)\leq p-\frac{2p}{p-1}<p-1=\ord(p),\]
so that the least-valuation term on the right-hand side of Equation (\ref{eqn2}) is $Av^p$. For the solvability it is necessary and sufficient that the right-hand side of Equation (\ref{eqn2}) is a square, by \cite[Lemma 2]{Fad61}, for which, it is necessary and sufficient that $v=Aw^2$, $w\in L$. Then the right-hand side of (\ref{eqn2}) takes the form
\[A^{p+1}w^{2p}-pB+A^2pw^2+\cdots=A^{p+1}w^{2p}(1-pw^{-2p}B_1+pw^{2-2p}B_2+\cdots),\]
with $p$-adic integral $B_1,B_2,\cdots$.  The square root of the right-hand side will be look like:
\[A^{(p+1)/2}w^p(1-pw^{-2p}B_1+pw^{2-2p}B_2+\cdots)^{1/2}=C_0w^p+C_1pw^{-p}+\gamma_1 p^2w^{-3p}+\gamma_2pw^{2-p},\]
where $C_0$ and $C_1$ are $p$-adic integral rational numbers, $\gamma_1$ and $\gamma_2$ are $p$-adic integral in $L$.

We estimate the orders of $\gamma_1p^2w^{-3p}$ and $\gamma_2pw^{2-p}$, taking into account that 
\[\ord(w)\leq \frac{1}{2}-\frac{1}{p-1}.\]
We have \[\ord(p^2w^{-3p})=2(p-1)-3p \ord(w)>\frac{p+1}{2},\]
\[\ord(pw^{2-p})=p-1-(p-2)\ord(w)>\frac{p+1}{2}.\]
So the equation for determining $z$ is 
\[z+a_2'z^2+\cdots=C_0w^p+C_1 pw^{-p}+\lambda^{(p+1)/2}\gamma,\]
for $\gamma\in L,\ord(\gamma)>0$. Hence,
\[z=\Phi(w^p)+C_1pw^{-p}+\lambda^{(p+1)/2}\gamma',\]
with $\gamma'\in L,\ord(\gamma')>0,$ $\Phi$-- a plynomial in $w^p$ with $p$-adic integral coefficients and with a zero constant term.

According to \cite[Lemmas 3, 4]{Fad61}, 
\[\Tr_{L/k}(\Phi(w^p))\equiv 0\mod \lambda^p,\quad \ord(\Tr_{L/k}(\gamma'))\geq 1,\]
and 
\begin{eqnarray*}
\ord(\Tr_{L/k}(pw^{-p}))&=&\ord(\Tr_{L/k}(-\lambda^{-1}(\lambda/w)^p))\\
&\geq& -1+\min(p\ord(\Tr_{L/k}(\lambda/w)), p-1+p\ord(\lambda/w))\\
&\geq &p-1.
\end{eqnarray*}
Therefore 
\[\Tr_{L/k}(z)\equiv\Tr_{L/k}(z^2)\equiv\cdots\equiv 0\mod \lambda^{(p+3)/2}.\]
Hence,
\[\N_{L/k}(x)=\N_{L/k}\left(\frac{r\delta}{r+1}+z\right)\equiv \left(\frac{r\delta}{r+1}\right)^{[L:k]}\mod \lambda^{(p+3)/2},\]
and consequently, $\N_{L/k}(x)$ falls into the subgroup generated by $u_i$, $i\geq (p+3)/2$.

\noindent{\bf V-B. } Suppose now $\ord(v)=1$, and $v=\lambda \alpha$ with $\ord(\alpha)=0$. In this case the right-hand side of Equation (\ref{eqn2}) is equal to 
\[A(p\lambda(\alpha-\alpha^p)+ \frac{p(p-1)}{2}\lambda^2\alpha^2+\cdots)+B\lambda^{p-1}.\]
We set 
\[A(p\lambda(\alpha-\alpha^p)+\frac{p(p-1)}{2}\lambda^2\alpha^2+\cdots)+B\lambda^{p-1}=u^2,\]
which must be satisfied in the field $L$ for the solvability of Equation (\ref{eqn2}). This case is broken up into three cases.

(a). $B\equiv 0\mod p$. Then
\[\ord(u)\geq p/2,\quad \ord(z)\geq p/2, \quad \ord(z^2)\geq p,\]
\[\ord(\Tr_{L/k}(z))\geq (p+1)/2,\quad \ord(\Tr_{L/k}(z^2))\geq p,\]
and
\[\N_{L/k}(x)=\N\left(\frac{r\delta}{r+1}+z\right)\equiv \left(\frac{r\delta}{r+1}\right)^{[L:k]}\cdot (1+\lambda^{(p+1)/2}b),\quad b\in k, \ord(b)\geq0.\]
Therefore, $\N_{L/k}(x)$ falls into the subgroup generated by $u_i$, $i\geq(p+1)/2$.

(b.)$\left(\frac{B}{p}\right)=+1$, so that $B=b^2,b\in \BQ_p$. In this case,
\[u=b\lambda^{(p-1)/2}\left(1-B^{-1}A\lambda(\alpha-\alpha^p)-B^{-1}A\frac{p-1}{2}\lambda^2\alpha^2-\cdots\right)^{1/2}=\]
\[b\lambda^{(p-1)/2}+c b\lambda^{(p+1)/2}(\alpha-\alpha^p)+\gamma\lambda^{(p+3)/2},\]
\[c\in \BQ_p,\quad \ord(c)=0, \quad \gamma\in L,\quad \ord(\gamma)\geq0.\]
According to \cite[Lemma 1]{Fad61},
\[z=b\lambda^{(p-1)/2}+c_1b\lambda^{(p+1)/2}(\alpha-\alpha^p)+\gamma'\lambda^{(p+3)/2}, \ord(\gamma')\geq0.\]
Note
\[\Tr_{L/k}(\alpha-\alpha^p)\equiv 0\mod \lambda.\]
Hence,
\[\N_{L/k}(x)=\N_{L/k}\left(\frac{r\delta}{r+1}+z\right)\equiv \left(\frac{r\delta}{r+1}\right)^{[L:k]}\cdot \left(1+[L:k]\left(\frac{r\delta}{r+1}\right)^{-1}b\lambda^{(p-1)/2}\right) \mod \lambda^{(p+3)/2}.\]
Thus, $\N_{L/k}(x)$ falls  into the subgroup generated by $u_{(p-1)/2}$ and $u_i$, $i\geq (p+3)/2$.

(c).  $\left(\frac{B}{p}\right)=-1$, so that $B=b^2$ for some $b$, belonging to the unramified quadratic extension $\BQ_{p^2}$ of $\BQ_p$. In this case,
as well as in the previous one,
\[z=b\lambda^{(p-1)/2}+c_1b\lambda^{(p+1)/2}(\alpha-\alpha^p)+\gamma'\lambda^{(p+3)/2}.\]
Now $\Tr(b\lambda^{(p-1)/2})=0$. But the trace of the second term can have order $(p+1)/2$,  for this it is sufficient to take, for example, $\alpha = b$, since then
\[b(\alpha-\alpha^p)=b^2(1-b^{p-1})=b^2(1-B^{(p-1)/2})\equiv 2B\mod (p).\]
Thus, in this case 
\[\N_{L/k}(x)= \left(\frac{r\delta}{r+1}\right)^{[L:k]}\N_{L/k}\left(1+\frac{r+1}{r\delta}z\right)\]
falls into the group spanned by $u_i,i\geq(p+1)/2$, and there exists an element with an non-trivial 
component for $u_{(p+1)/2}$.

\noindent{\bf V-C.}  Suppose now that $\ord(v)> 1, v = \lambda w.$ In this case, the right-hand side of
Equation (\ref{eqn2}) has the form:
\[Ap \lambda w+Cp \lambda^2 w^2+\cdots +B\lambda^{p-1}=u^2.\]
Again, we shall consider three subcases.

(a). $ B \equiv  0 \mod (p)$.
In this case 
\[\ord(u)> p/2, \quad \ord(z)>p/2,\quad \ord(\Tr(z))\geq(p+1)/2,\]
 so that
\[\N_{L/k}(x)=\N_{L/k}\left(\frac{r\delta}{r+1}+z\right)\]
 falls into the  subgroup spanned by  $u_i, i\geq (p+1)/2.$ For any $i$ in the interval   $(p+1)/2\leq i\leq p-2$, one can take $L=k,w=-A\lambda^{2i-p}$.
With this choice,
\[u^2=A^2\lambda^{2i}+\cdots, \quad u\in k,\  \ord(u)=\ord(z)=a,\]
\[x=\N_{L/k}(x)=\frac{r\delta}{r+1}(1+c\lambda ^i),\quad  c\not\equiv 0\mod(\lambda). \]
Thus, in the group generated by $\N_{L/k}(x)$, there are elements with non-zero exponents for $u_i$ for each $i$ in the interval $(p+1)/2 \leq i\leq p-2 $.

(b). $\left(\frac{B}{p}\right)=+1,$ and so $ B=b^2,b\in \BQ_p$. In this case,
\[u^2=b^2\lambda^{p-1}+Ap\lambda w+\cdots;\]
\[u=b\lambda^{(p-1)/2}-\frac{1}{2}Ab^{-1}\lambda^{(p+1)/2}w+\gamma \lambda^{(p+3)/2}w;\]
\[\gamma\in L,\  \ord(\gamma)>0.\]
Then
\[z=b\lambda^{(p-1)/2}-\frac{1}{2}Ab^{-1}\lambda^{(p+1)/2}w+\gamma' \lambda^{(p+3)/2}w,\quad \gamma'\in L,\ \ord(\gamma')\geq 0,\]
\[\Tr_{L/k}(z)=[L:k]b\lambda^{(p-1)/2}+c,c\in k,\ord(c)\geq (p+3)/2,\]
so $\N_{L/k}(x)$ falls into the subgroup spanned by $u_{(p-1)/2}$ and $u_i, i\geq (p+3)/2$.

For any $i$ from the interval  $(p+3)/2\leq i\leq p-2$,  we can take $ L=k, w=\lambda^{i-(p+1)/2}$. With this choice 
\[z =\Tr_{L/k}(z)=b\lambda^{(p-1)/2}-\frac{1}{2}Ab^{-1}\lambda^{i}+c,\quad c\in k,\ \ord(c)>i,\]
and the exponents for $u_{(p-1)/2}$ and $u_i$ in $\N_{L/k}(x)$ turn out to be different from $0.$

(c).   $\left(\frac{B}{p}\right)=-1$, so that $B=b^2$ for some $b$, belonging to the unramified quadratic extension $\BQ_{p^2}$ of $\BQ_p$.  As in the previous case,
\[z=b\lambda^{(p-1)/2}-\frac{1}{2}Ab^{-1}\lambda^{(p+1)/2}w+\gamma' \lambda^{(p+3)/2}w,\quad \ord(\gamma')>0,\]
Now we have $L=k\cdot\BQ_{p^2}$ and
\[\Tr_{L/k}(z)=-\frac{1}{2}A\lambda^{(p+1)/2}\Tr_{L/k}(b^{-1}w)+\lambda^{(p+3)/2}\Tr_{L/k} (\gamma' w)\equiv 0\mod \lambda^{(p+3)/2},\]
and  
\[\N_{L/k}(x)= \left(\frac{r\delta}{r+1}\right)^{[L:k]}\N_{L/k}\left(1+\frac{r+1}{r\delta}z\right)\]
falls in a group spanned by $u_i,i\geq (p+3)/2$.  We can take, for any $i$ from the interval $(p+3)/2\leq i\leq p-2$, that  $L=k\cdot \BQ_{p^2}, [L:k]=2, $ and $w=b\lambda^{i-(p+1)/2}.$ Then 
\[\Tr_{L/k}(z)=-A\lambda^i+c',\quad \ord(c')>i,\]
and $\N_{L/k}(x)$ has a non-zero exponent for $u_i$.

To summarize, in all the cases examined,
$\N_{L/k}(x)$ falls into the subgroup of the group $k^{\times}/k^{\times p} $ spanned by $u_i, (p-1)/2\leq    i   \leq p.$ This subgroup is invariant under the Galois group $\Gal(K/\BQ)$, and has multiplicity one for each occurring eigen-space under the Galois action. Analysis of all cases gives the following therorem.

Recall we have set
\[B=p^{-1}\left(\left(\frac{r^r\delta^{r+1}}{(r+1)^{r+1}}\right)^{p-1}-1\right)\cdot \frac{2r\delta^2}{(r+1)^3}.\]
\begin{thm}\label{desc-at-p1}
Suppose the triple $(r,s,t)=(r,1,p-r-1)$ with $1\leq r\leq p-2$ and $\delta\in \BQ^\times/\BQ^{\times p}$ is an integer such that $p\nmid \delta$. The Kummer image $\Im(\delta_\Pi)$ is the subgroup of $K_\Pi^\times/K_\Pi^{\times p}$ generated by $u_i$, $(p+3)/2\leq i\leq p$, and, in addition to them, $u_{(p-1)/2}$ if $\left(\frac{B}{p}\right)=+1$, and $u_{(p+1)/2}$ if $\left(\frac{B}{p}\right)=-1$ or $B\equiv 0\mod p$.
\end{thm}
\begin{proof}
The theorem is a summary of all the above examined cases. 
\end{proof}

\subsection{Selmer groups for the Fermat Jacobian varieties}\label{Selmer-section}
Denote $S=\{V : V\mid p\delta\}$, the set of places of $K$ dividing $p\delta$, and  $R=\CO_{K}[\frac{1}{p\delta}]$. Then $R^\times$ is the group of $S$-units in $K$. The torsion subgroup of $R^\times$ is the group of roots of unity, and by Dirichlet theorem, the free part of $R^\times$ has $\BZ$-rank $|S|+(p-3)/2$.

Taking cohomology from the short exact sequence
\[1\ra \mu_p\ra \BG_m\ra \BG_m\ra 1,\]
we have the exact sequence
\begin{equation}\label{class-group-vs-ST-group}
1\ra R^\times/R^{\times p}\ra \RH_\et^1(R,\mu_p)\ra \RH^1_\et(R,\BG_m)[p]=\Cl(R)[p].
\end{equation}
We identify $J_\delta[\Pi]$ with $\mu_p$ as in (\ref{id}) through the following composition
\[J_\delta[\Pi]\xrightarrow{e_\Pi(D_0,\cdot)}\mu_p,\]
and then the \'etale cohomology group 
\begin{equation}\label{ident}
\RH^1_\et(R,J_\delta[\Pi])\simeq\{x\in K^\times/K^{\times p}: \ord_V(x)\equiv 0\mod p \textrm{ for all places }V\nmid p\delta\}. 
\end{equation}
Since each class of the  Selmer group $\Sel({J_{\delta}}_{/K},\Pi)\subset \RH^1(K, J_\delta[\Pi])$ is unramified at any place $V\nmid p\delta$, the  Selmer group $\Sel({J_{\delta}}_{/K},\Pi)$ is a subgroup of $\RH^1_\et(R,J_\delta[\Pi])$. 

Let $\alpha:\Sel\left({J_\delta}_{/K}, \Pi\right)\ra \Cl(R)[p]$ be the composite of 
\[\Sel\left({J_\delta}_{/K}, \Pi\right)\hookrightarrow \RH^1_\et(R,J_{\delta}[\Pi])\ra \Cl(R)[p].\]
Then
\[\Ker(\alpha)=\{x\in R^\times/R^{\times p}: x_V\in \Im(\delta_V) \textrm{ for all } V\}.\]
and we have an exact sequence
\begin{equation}\label{exact1}
0\ra \Ker(\alpha)\ra \Sel\left({J_\delta}_{/K}, \Pi\right)\xrightarrow{\alpha} \Cl(R)[p].
\end{equation}

\subsubsection{General bounds}

Let $\sigma:(\BZ/p\BZ)^\times\ra \Gal(K/\BQ)$ be the Artin reciprocity map such that $\omega^{\sigma_a}=\omega^a$ for $a\in (\BZ/p\BZ)^\times$. Let $g\in (\BZ/p\BZ)^\times$ be a primitive root modulo $p$. For $2\leq i\leq p-3$ even, define the cyclotomic units
\[E_i=\prod_{a\in (\BZ/p\BZ)^\times} \left(\omega^{\frac{1-g}{2}}\frac{1-\omega^g}{1-\omega}\right)^{a^i\sigma_a^{-1}}\in \CO_F^\times/\CO_F^{\times p}(i).\]

\begin{lem}\label{p-power}
An element in $\CO_{K,\Pi}^\times$ is a $p$-th power if and only if it is congruent to rationals modulo $\Pi^{p+1}$, i.e. $\CO_{K,\Pi}^{\times p}=\BF_p^\times(1+\Pi^{p+1}\CO_{K,\Pi})$.
\end{lem}
\begin{proof}
First we show $\BF_p^\times(1+\Pi^{p+1}\CO_{K,\Pi})\subset \CO_{K,\Pi}^{\times p}$. Since $\BF_p^\times$ has order $p-1$, any element in $\BF_p^\times$ is a $p$-th power. If $\ord(x)>p$, then the formal power series
\[(1+x)^{1/p}=\sum_{n=0}^\infty \left(\begin{matrix}1/p\\n\end{matrix}\right)x^n\]
coverges, and hence $1+x$ is a $p$-th power.  The lemma follows by noting the equality of the indices
\[[\CO_{K,\Pi}^\times: \CO_{K,\Pi}^{\times p}]=[\CO_{K,\Pi}^\times: \BF_p^\times(1+\Pi^{p+1}\CO_{K,\Pi})]=p^p.\]
\end{proof}

For even $i$ with $2\leq i\leq p-3$, let $B_i$ be the $i$-th Bernoulli number.
\begin{lem}\label{non-divisibility}
If $\ord_p(B_i)=0$, then $E_i$ is not a $p$-th power in $K_{\Pi}^\times$. 
\end{lem}
\begin{proof}
Under the hypothesis of the lemma, we have
\[\ord_p(L_p(1,\mu^i))=\ord_p(B_i)=0.\]
By \cite[Proposition 8.12]{Washington-gtm83},
\[\ord_p(\log_p E_i)=\frac{i}{p-1}.\]
Then $E_i$ has $p$-adic expansion
\[E_i=a+b\Pi^i+c\Pi^{i+2}+\cdots\]
where the coefficients $a,b,c,\cdots$ are rational $p$-adic units. Since $i\leq p-3$, by Lemma \ref{p-power}, $E_i$ is not a $p$-th power in $K_\Pi^\times$.
\end{proof}
Denote $i(p)$ the index of irregularity of $p$, i.e. the number of Bernoulli numbers $B_{i}$, with $i$ even and $2\leq i\leq p-3$, such that $\ord_p( B_{i})>0$, and $k(\delta)$ the number of distinct prime ideals of $K$ dividing $\delta$. 
\begin{thm}\label{selmer-bound}
Suppose $p\geq 5$ and suppose the triple $(r,s,t)=(r,1,p-r-1)$ with $1\leq r\leq p-2$ and $\delta$ is a $p$-th power-free integer such that $p\nmid \delta$. Then
\[\dim_{\BF_p}\Sel\left({J_\delta}_{/K}, \Pi\right)\leq k(\delta)+\frac{p-3}{4}+i(p)+\dim_{\BF_p}\Cl(K)[p].\]
\end{thm}
\begin{proof}
It follows from the exact sequence (\ref{exact1}) that 
\[\dim_{\BF_p}\Sel\left({J_\delta}_{/K}, \Pi\right)\leq \dim_{\BF_p}\ker(\alpha)+\dim_{\BF_p}\Cl(K)[p].\]
Here we note that $\Cl(R)$ is a natural quotient of $\Cl(K)$.
The dimension of $R^\times/R^{\times p}$ is $|S|+(p-1)/2=k(\delta)+(p+1)/2$. We shall find out certain nontrivial classes $x\in (R^\times/R^{\times p})(i)$ such that $x\not \in \Ker(\alpha)$ for $0\leq i\leq (p-3)/2$.  For $i=0,1$, by Theorem \ref{desc-at-p1}, we may take $x=p, \omega$ respectively. If $i$ is even and $2\leq i\leq (p-3)/2$ and if $p\nmid B_i$, then by Lemma \ref{non-divisibility}, $E_i$ is not a $p$-th power in $K_\Pi^\times$ and generates the $i$-th eigenspace $U_i$, and again by Theorem \ref{desc-at-p1} we may take $x=E_i$. Then
\begin{eqnarray*}
\dim_{\BF_p}\Ker(\alpha)&\leq& k(\delta)+\frac{p-3}{2}-\#\{k: 2\leq 2k\leq (p-3)/2, p\nmid B_{2k}\}\\
&\leq &k(\delta)+\frac{p-3}{2}-\left[\frac{p-3}{4}\right]+i(p)\\
&\leq &k(\delta)+\frac{p-3}{4}+i(p)
\end{eqnarray*}
\end{proof}
By Ribet's converse to Herbrand's theorem \cite{Ribet76},
\[i(p)\leq \dim_{\BF_p}\Cl(K)[p].\]
Therefore, we have
\begin{coro}
Let the notations be as in Theorem \ref{selmer-bound}.  Then
\[\dim_{\BF_p}\Sel\left({J_\delta}_{/K}, \Pi\right)\leq k(\delta)+\frac{p-3}{4}+2\dim_{\BF_p}\Cl(K)[p].\]
\end{coro}

\subsubsection{Structures of Selmer groups under regular conditions}
In this section we assume $p\geq 5$ is regular, i.e. $\Cl(K)[p]=0$, and $\delta$ is a $p$-th power-free integer satisfying the Hypothesis $(\mathrm{\bf{I}})$. Since $\Cl(R)$ is a natural quotient of $\Cl(K)$, we have $\Cl(R)[p]=0$. Then by the exact sequence (\ref{exact1}), 
\[\Sel\left({J_\delta}_{/K}, \Pi\right)=\Ker(\alpha)=\{x\in R^\times/R^{\times p}: x_V\in \Im(\delta_V) \textrm{ for all } V\}.\]
We first determine the $\BZ_p[\Gal(K/\BQ)]$-module structure of $R^\times/R^{\times p}$ and then use this to determine the structure of the Selmer group $\Sel\left({J_\delta}_{/K}, \Pi\right)$.

Recall $F$ is the maximal real subfield of $K$. Since $p$ is regular, $\Cl(F)[p]=0$ and denote the class number $h=\# \Cl(F)$.
For any place $V$ of $K$ above $\delta$, denote $\fP_V$ the corresponding prime ideal.  Let $\loc_\Pi: K^\times \ra K_\Pi^\times/K_\Pi^{\times p}$ be the localization map at $\Pi$.
\begin{lem}\label{local-pi}
For any place $V\mid \delta$, the ideal $\fP_V^h$ is principal and there exists a generator $\varpi_V\in K$ satisfying 
\begin{itemize}
\item[1.] $\varpi_V\in \CO_F$;
\item[2.] The localization $\loc_\Pi(\varpi_V)\in U_{p-1}$.
\end{itemize}
\end{lem}
\begin{proof}
Let $v$ be place of $F$ restricted from $V$ and let $\fP_v$ be the corresponding prime ideal. Then $\fP_v^h$ is a principal ideal generated by an element $\varpi_V\in \CO_F$. Under the Hypothesis $(\mathrm{\bf{I}})$, $V$ is inert in $K/F$, and hence $\fP_V^h$ is also principal and generated by $\varpi_V$. Since $p$ is regular,  by \cite[Theorem 5.34]{Washington-gtm83},  $\ord_p(B_i)=0$ for all even $i$ with $2\leq i\leq p-3$. Then by Lemma \ref{non-divisibility}, each localization $\loc_\Pi(E_i)$ is nontrivial and generates $U_i$.
Since $\varpi_V$ is real, $\loc_\Pi(\varpi_V)$ lies in the direct sum of even eigen-spaces. Therefore by suitably multiplying by the units $E_i$ with $i$ even and $2\leq i\leq p-3$, we may make  a choice of the generator $\varpi_V$ so that $\loc_\Pi(\varpi_V)\in K_\Pi^\times/K_\Pi^{\times p}(0)=U_0\times U_{p-1}$. Since $\varpi_V$ is a $p$-adic unit, we conclude that $\loc_\Pi(\varpi_V)\in U_{p-1}$.
\end{proof}

For $V=\Pi$, denote $\varpi_\Pi=p$. We note for each $V\in S$, $\varpi_V$ is not a $p$-th power  in $K^\times$ since $p\nmid h(p-1)$ and $\varpi_V$ is a $V$-unit with positive valuation at $V$.

\begin{prop}\label{units}
The $\BF_p$-space $R^\times/R^{\times p}$ has a basis consisting of
\[ \varpi_V \text{ with }V\in S, \quad \omega \text{ and }E_i \text{ with $i$ even and $2\leq i\leq p-3$}.\]
\end{prop}
\begin{proof}
First we note, by the Dirichlet unit theorem,  $R^\times/R^{\times p}$ has $\BF_p$-dimension $|S|+(p-1)/2$. Since $\Cl(F)[p]=0$, by \cite[Theorem 8.14, Corollary 8.15]{Washington-gtm83}, $E_2,\cdots, E_{p-3}$ form a basis of $\CO_F^\times/\CO_F^{\times p}$.  Since $\CO_K^\times=W\CO_F^\times$, where $W$ is the subgroup of roots of unity, we have 
\[\CO_K^\times/\CO_K^{\times p}(1)=\langle\omega\rangle,\quad \CO_K^\times/\CO_K^{\times p}(i)=\langle E_i\rangle,\]
for even $i$ with $2\leq i\leq p-3$.
On the other hand, $\varpi_V, V\in S$, are linearly independent $S$-units which are also linearly independent with the units.
The proposition follows by counting the cardinality.
\end{proof}

\begin{lem}\label{p-power2}
Suppose $V$ is a place of $K$ that is inert in $K/F$. If $x\in F^\times$ is a $V$-adic unit, then $x\in K_V^{\times p}$.
\end{lem}
\begin{proof}
First we note necessarily $V\nmid p$. Let $\ell$ be the rational prime under $V$ and let $f$ be the order of $\ell$ in $\BF_p^\times$. Then $f$ is the inertia degree $[\BF_V:\BF_p]$. Let $v$ be the place of $F$ restricted from $V$ and let $\varpi$ be a uniformizer of the local field $F_v$. We have $\CO_{F,v}^\times=\BF_v^\times(1+\varpi\CO_{F,v})$. Since $V$ is inert in $K/F$, $\BF_v^\times$ has cardinality $\ell^{f/2}-1$ which is prime to $p$. Then taking $p$-th power induces isomorphisms on both $\BF_v^\times$ and $1+\varpi\CO_{F,v}$ and hence $\CO_{F,v}^\times=\CO_{F,v}^{\times p}$. In particular, as an element of $\CO_{F,v}^\times$, $x\in K_V^{\times p}$.
\end{proof}

\begin{thm}\label{Selmer}
Suppose $p\geq 5$ is a regular prime, the triple $(r,s,t)=(r,1,p-r-1)$ with $1\leq r\leq p-2$ and $\delta$ is a $p$-th power-free integer satisfying the Hypothesis $(\mathrm{\bf{I}})$. If $p\equiv 1\mod 4$, the Selmer group  $\Sel\left({J_\delta}_{/K}, \Pi\right)$ has a basis consisting of
\[\varpi_V \text{ with }V\mid \delta \text{ and } E_i \text{ with $i$ even and }  (p+3)/2\leq i\leq p-3,\]
and in addition to them, $E_{(p-1)/2}$ if $\left(\frac{B}{p}\right)=+1$. If $p\equiv 3\mod 4$, the Selmer group  $\Sel\left({J_\delta}_{/K}, \Pi\right)$ has a basis consisting of
\[\varpi_V\text{ with }V\mid \delta \text{ and } E_i \text{ with $i$ even and } (p+5)/2\leq i\leq p-3,\]
and in addition to them, $E_{(p+1)/2}$ if $\left(\frac{B}{p}\right)=-1$ or $B\equiv 0\mod p$. Here we accept the convention that $\left(\frac{0}{p}\right)=-1$ and 
\[B=p^{-1}\left(\left(\frac{r^r\delta^{r+1}}{(r+1)^{r+1}}\right)^{p-1}-1\right)\cdot \frac{2r\delta^2}{(r+1)^3}.\]
\end{thm}
\begin{proof}
The generators of $R^\times/R^{\times p}$ are determined by Proposition \ref{units}. We shall check the local Kummer conditions of Theorems \ref{desc-off-p} and \ref{desc-at-p1} for the generators of $R^\times/R^{\times p}$. 

First note $p\geq 5$ and $\omega$ generates $U_1$. It follows immediately from Theorem \ref{desc-at-p1} that $\omega\not\in \Im(\delta_\Pi)$. Next we consider the generators $\varpi_V, V\in S$. If $V=\Pi$, then $\varpi_\Pi=p$ generates $U_0$ and hence $\varpi_\Pi\not\in \Im(\delta_\Pi)$ by Theorem \ref{desc-at-p1}. 

Suppose $V\mid \delta$ and we consider the generator $\varpi_V$. At the place $\Pi$, by Lemma \ref{local-pi},  $\varpi_V\in U_{p-1} \subset \Im(\delta_\Pi)$. The last inclusion is ganranteed by $p\geq 5$ and Theorem \ref{desc-at-p1}. 

 Assume $\ell$ is the rational prime under $V$ and $\ord_\ell(\delta)=m$ with $1\leq m\leq p-1$. If $1=px+hy$ for some integers $x,y$, then 
 \[\delta\equiv \varpi_V^{my} \delta_0\mod K_V^{\times p},\]
 where $\delta_0\in F^\times$ is a $V$-adic unit. Then by Lemma \ref{p-power2}, $\delta_0\in K_V^{\times p}$ and hence by Theorem \ref{desc-off-p},
 \[\varpi_V\in \langle\varpi_V\rangle=\langle \delta\rangle=\Im(\delta_V).\]

At the places $W\neq \Pi \text{ or }V$, $\varpi_V$ is a $W$-adic unit and hence 
\[\varpi_V\in \RH^1_\ur(K_W,J_\delta[\Pi])\simeq\BF_W^\times/\BF_W^{\times p}.\]
Moreover, if $W\mid \delta$, by Lemma \ref{p-power2}, $\varpi_V\in K_W^{\times p}$. It follows from Theorem \ref{desc-off-p} that $\varpi_V\in \Im(\delta_{W}) $. 

Thus we see the generators $\omega, \varpi_\Pi$ don't lie in the Selmer group $\Sel\left({J_\delta}_{/K}, \Pi\right)$ and the generators $\varpi_V, V\mid \delta,$ all lie in the Selmer group $\Sel\left({J_\delta}_{/K}, \Pi\right)$.

Finally we consider the generators $E_2,E_4,\cdots,E_{p-3}$. For all $V\nmid p$, clearly for all even $2\leq i\leq p-3$, as units, 
\[E_i\in \RH^1_\ur(K_V,J_\delta[\Pi])\simeq \BF_V^\times/\BF_V^{\times p}.\]
Moreover, if $V\mid \delta$, by Lemma \ref{p-power2}, $E_i\in K_V^{\times p}$. It follows from Theorem \ref{desc-off-p} that $E_i\in \Im(\delta_{V}) $. 

At the place $\Pi$, by Lemma \ref{non-divisibility}, the localization of each $E_i$ is nontrivial and generates the $i$-th eigen-space $U_i$. It follows from Theorem \ref{desc-at-p1} that the generators $E_i\in \Im(\delta_\Pi)$ are exactly those listed in the Theorem.
Thus we conclude the theorem.

\end{proof}

Recall $k(\delta)$ denotes the number of distinct prime ideals of $K$ dividing $\delta$. As an immediate consequence we have the following dimension formulae.
\begin{thm}\label{Sel-dim}
Let the notations be as in Theorem \ref{Selmer}. We have
\[\dim_{\BF_p}\Sel\left({J_\delta}_{/K}, \Pi\right)=k(\delta)+\left\{\begin{aligned}
\frac{1}{4}\left(p-3+2\left(\frac{B}{p}\right)\right),\quad \textrm{ if } p\equiv 1\mod 4;\\
\frac{1}{4}\left(p-5-2\left(\frac{B}{p}\right)\right),\quad \textrm{ if } p\equiv 3\mod 4.
\end{aligned}\right.\]

\end{thm}

\begin{remark}
If $p=3$, then $r=s=t=1$ and both $F_\delta: x^3+y^3=\delta$ and $C_\delta:y^3=x(\delta-x)$ are elliptic curves over $\BQ$. The $3$-Selmer groups of these elliptic curves via $3$-descent method have been extensively investigated by Satg\'e \cite{Satge}. For the Jacobian varieties of the Fermat curves $X^p+Y^p=1$,  Faddeev \cite{Fad61} calculated the Selmer groups of these Jacobian varieties when $p\geq 5$ is regular.
\end{remark}
For fixed $1\leq r\leq p-2$, denote
\[S_\delta=S_{r,1,p-r-1;\delta}=\dim_K {J_\delta(K)\otimes_{\CO_K}K}+\rk_{\BF_p}\Sha\left({J_\delta}_{/K}\right)[\Pi].\]
We have the following
\begin{coro}\label{selmer-rank}
Let the notations be as in Theorem \ref{Selmer}. We have
\[S_\delta=k(\delta)+\left\{\begin{aligned}
\frac{1}{4}\left(p-7+2\left(\frac{B}{p}\right)\right),\quad \textrm{ if } p\equiv 1\mod 4;\\
\frac{1}{4}\left(p-9-2\left(\frac{B}{p}\right)\right),\quad \textrm{ if } p\equiv 3\mod 4.
\end{aligned}\right.\]
\end{coro}
\begin{remark}
This is compatible with \cite[Proposition 4.1]{Gro-Roh78} for the case $\delta=1$. Note there is a sign error in the second formula of \cite[Proposition 4.1]{Gro-Roh78}.
\end{remark}

\begin{proof}
We have the exact sequence
\[0\ra J_\delta(K)/[\Pi] J_\delta(K)\ra \Sel({J_\delta}_{/K},\Pi)\ra\Sha\left({J_\delta}_{/K}\right)[\Pi]\ra 0,\]
and
\[\dim_{\BF_p}J_\delta(K)/\Pi J_\delta(K)=\dim_{\BF_p} J_\delta(K)_\tor/[\Pi] J_\delta(K)_\tor+\dim_{\BF_p}  {J_\delta(K)}_\mathrm{fr} /[\Pi]{J_\delta(K)}_\mathrm{fr},\]
where  $ {J_\delta(K)}_\tor$ and $ {J_\delta(K)}_\mathrm{fr}$ denotes the torsion part and the free part respectively.
The torsion group $J_\delta(K)[\Pi]$ is nonzero and generated by the class $D_0$ (for example see \cite{Greenberg80}). Suppose the $\Pi$-primary compoent $J_\delta(K)[\Pi^\infty]=J_\delta[\Pi^n]$ for some $n\geq 1$. Since the abelian varieties $J_\delta$ have CM by $\CO_K$, $J_\delta[\Pi^n]\simeq \CO_K/(\Pi^n)$. Then
\[J_\delta(K)_\tor/\Pi J_\delta(K)_\tor\simeq J_\delta[\Pi^n]/\Pi J_\delta[\Pi^n]\simeq \CO_K/(\Pi)=\BF_p.\]
Since
\[\dim_{\BF_p}  {J_\delta(K)}_\mathrm{fr} /\Pi{J_\delta(K)}_\mathrm{fr}=\dim_K J_{\delta}(K)\otimes_{\CO_K}K,\]
we have
\[S_\delta= \dim_{\BF_p}\Sel({J_\delta}_{/K},\Pi)-1,\]
and the corollary follows.
\end{proof}

\section{The Parity Conjecture}
We are ready to verify the parity conjecture (\ref{pc}) for the Jacobian varieties of Fermat curves under the regular condition and the Hypothesis ($\mathrm{\bf{I}}$).

Denote
\[S_\pm=\left\{ \text{rational primes } \ell \mid \delta: \left(\frac{\ell}{p}\right)=\pm 1\right\}.\]
\begin{lem}\label{mod2}
Let $\delta$ be a $p$-th power-free integer satisfying the Hypothesis $(\mathrm{\bf{I}})$. Then
\[k(\delta)\equiv \# S_-\mod 2.\]
\end{lem}
\begin{proof}
Let $\ell\mid \delta$ be a rational prime and $f_\ell$ the order of $\ell \mod p$. Since $\ell$ is unramified in $K/\BQ$ and $f_\ell$ is the inertia degree of $\ell$ in $K/\BQ$, there are $(p-1)/f_\ell$ places of $K$ above  $\ell$. If $\ell\in S_+$ resp. $\ell\in S_-$, then $(p-1)/f_\ell$ is even resp. odd. Then
\[k(\delta)=\sum_{\ell\mid \delta} (p-1)/f_\ell \equiv \# S_-\mod 2.\]

\end{proof}

\begin{thm}\label{parity}
Suppose $p\geq 5$ is a regular prime and $(r,s,t)$ is a triple with $r,s,t>0$ and $r+s+t=p$. Suppose $\delta$ is a $p$-th power-free integer satisfying the Hypothesis $(\mathrm{\bf{I}})$. Then \[\epsilon_{r,s,t;\delta}=(-1)^{S_{r,s,t;\delta}}.\]
\end{thm}
\begin{proof}
As is noted that if 
\[(r',s',t')=h(r,s,t)+p(i,j,k),\]
with $r',s',t'>0,$ and $r'+s'+t'=p$. Then we have a commutative diagram
\[\xymatrix{&F_\delta\ar[ld]_{\varphi_{r,s,t}}\ar[rd]^{\varphi_{r',s',t'}}&\\C_{r,s,t;\delta}\ar[rr]&&C_{r',s',t';\delta}\ ,}\]
where the horizontal arrow is
\[(x,y)\mapsto(x,y^hx^i(\delta-x)^j).\]
Since $h$ is prime to $p$, this map is a birational equivalence. 
Without loss of generality, we may suppose $s=1$ and take $1\leq r\leq p-2$.

Let $d\in \BZ/p\BZ$ be such that
\[d=p^{-1}\left(\left(\frac{r^r\delta^{r+1}(r+1)^p}{{(r+1)^{r+1}}}\right)^{p-1}-1\right)\mod p.\]
Since $(\BZ/p^2\BZ)^\times$ has cardinality $p(p-1)$, 
\[(r+1)^{p(p-1)}\equiv 1\mod p^2,\]
and hence
\begin{eqnarray*}
B&=&p^{-1}\left(\left(\frac{r^r\delta^{r+1}}{(r+1)^{r+1}}\right)^{p-1}-1\right)\cdot \frac{2r\delta^2}{(r+1)^3}\\
&\equiv& d\frac{2r(r+1)\delta^2}{(r+1)^4}\mod p.
\end{eqnarray*}
On the other hand, put
\[s=\left\{\begin{aligned}
\frac{1}{4}\left(p-7+2\left(\frac{B}{p}\right)\right),\quad \textrm{ if } p\equiv 1\mod 4;\\
\frac{1}{4}\left(p-9-2\left(\frac{B}{p}\right)\right),\quad \textrm{ if } p\equiv 3\mod 4.
\end{aligned}\right.\]
Here we accept the convention $\left(\frac{0}{p}\right)=-1$ as before.
It is straight-forward to verify that
\begin{equation}\label{link}
(-1)^s=-\left(\frac{B}{p}\right)\left(\frac{2}{p}\right).
\end{equation}

Since $p\nmid \delta$, in the notation as in \cite[Theorem 1.1]{Shu21jnt}, we have
\[u(r^r(-r-1)^{p-(r+1)\delta^{r+1}})\geq 1.\]
First suppose $u(r^r(-r-1)^{p-(r+1)\delta^{r+1}})\geq 2$, or equivalently
\[B\equiv d \equiv 0\mod p.\]
By \cite[Theorem 1.1]{Shu21jnt}, Corollary \ref{selmer-rank}, Lemma \ref{mod2} and (\ref{link}), we have
\[\epsilon_{r,1,p-r-1;\delta}=(-1)^{\# S_-}\left(\frac{2}{p}\right)=(-1)^{k(\delta)}\left(\frac{2}{p}\right)=(-1)^{S_\delta}.\]
Next we suppose $u(r^r(-r-1)^{p-(r+1)\delta^{r+1}})=1$. Then
\[B\equiv d\frac{2r(r+1)\delta^2}{(r+1)^4}\neq 0\mod p.\]
Again by  \cite[Theorem 1.1]{Shu21jnt}, Corollary \ref{selmer-rank}, Lemma \ref{mod2} and (\ref{link}), we have
\begin{eqnarray*}
\epsilon_{r,1,p-r-1;\delta}&=&(-1)^{\# S_-+1}\left(\frac{-r(p-(r+1)d)}{p}\right)\\
&=&(-1)^{k(\delta)+1}\left(\frac{-r(p-(r+1)d)}{p}\right)=(-1)^{k(\delta)+1}\left(\frac{2B}{p}\right)=(-1)^{S_\delta}.
\end{eqnarray*}
This completes the proof of the theorem.
\end{proof}

For a rational prime $p$ and real $x>0$, denote 
\[\Sigma(p,x)=\{\text{rational primes $\ell\leq x$ satisfying the Hypothesis ($\mathrm{\bf{I}}$)}\}.\]
The set $\Sigma(p,x)$ is Frobenian in the sense of \cite[(1.4)]{Serre74}.
\begin{prop}\label{density}
Among the set of all rational primes, the subset $\Sigma(p,x)$ has natrual density $1-2^{-n}$ where $n=\ord_2(p-1)\geq 1$.
\end{prop}
\begin{proof}
Let $\ell\neq p$ be a rational prime and let $f$ be the order of $\ell \mod p$.  It follows from \cite[Lemma 3.3]{Shu21jnt} that the rational primes $\ell$ satisfying the Hypothesis ($\mathrm{\bf{I}}$) are exactly those with $f$ even, or equivalently, the Frobenius element $\Fr_\ell$ has even order in $\Gal(K/\BQ)$. Since $\Gal(K/\BQ)\simeq (\BZ/p\BZ)^\times$, it is easy to see that the fraction of Galois automorphisms with even order is $1-2^{-n}$ where $n=\ord_2(p-1)$. Then the proposition follows by the Chebotarev's density theorem.
\end{proof}
Fix the triple $(r,s,t)$ with $r,s,t>0$ and $r+s+t=p$ and define
\[N_{r,s,t;\delta}(p,X)=\#\{\delta < X : |\delta|\leq X, \epsilon_{r,s,t;\delta}=S_{r,s,t;\delta}\}.\]

\begin{thm}
Suppose $p\geq 5$ is a regular prime and $(r,s,t)$ is a triple with $r,s,t>0$ and $r+s+t=p$. As $X\ra \infty$, we have
\[N_{r,s,t;\delta}(p,X)\gg \frac{X}{\log^{2^{-n}}X},\]
where $n=\ord_2(p-1)$.
\end{thm}
\begin{proof}
We just count the number of integers $\delta$ appearing in Theorem \ref{parity-intro}, say
\[S(X)=\sum_{|\delta|\leq X}1\]
where $\delta$ runs through $p$-th power-free integers satisfying the Hypothesis $(\mathrm{\bf{I}})$. In view of Proposition \ref{density}, the rational primes satisfying the Hypothesis $(\mathrm{\bf{I}})$ is Frobenian of natural density $1-2^{-n}$. Applying an variant  of Ikehara's Tauberian theorem \cite{Delang56} \cite[Theorem 2.4]{Serre74}, it follows that
\[S(X)\sim c\frac{X}{\log^{1-(1-2^{-n})}X}\sim c\frac{X}{\log^{2^{-n}}X},\]
where $c$ is a constant. The theorem follows immediately since $N_{r,s,t;\delta}(p,X)\geq S(X)$.
\end{proof}


\end{document}